\documentclass[a4paper,12pt]{article}
\usepackage{graphicx}
\usepackage{amsmath,amssymb,amsthm}
\usepackage{cases}
\usepackage{amsfonts}
\allowdisplaybreaks[4]
\usepackage{amscd}
\usepackage[all]{xy}
\usepackage{color}
\usepackage{comment}
\usepackage{url}

\DeclareMathOperator{\Ker}{Ker}
\DeclareMathOperator{\Hom}{Hom}

\DeclareMathOperator{\Homo}{H}

\DeclareMathOperator{\sol}{Sol}

\DeclareMathOperator{\cone}{cone}
\newcommand{\barsigma}{\overline{\sigma}}

\newcommand{\R}{\mathbb{R}}

\newcommand{\C}{\mathbb{C}}
\newcommand{\DD}{\mathcal{D}}

\newcommand{\Q}{\mathbb{Q}}
\newcommand{\Z}{\mathbb{Z}}

\newcommand{\ii}{\sqrt{-1}}
\newcommand{\s}{\sigma}

\newcommand{\bs}{\barsigma}

\newcommand{\HdR}{\mathcal{H}_{dR}}
\newcommand{\Db}{\mathfrak{D}\mathfrak{b}}

\def\pd#1{\partial_{#1}}

\theoremstyle{plain}
\newtheorem{theorem}{Theorem}[section]

\newtheorem{prop}[theorem]{Proposition}

\newtheorem{example}[theorem]{Example}

\newtheorem{rem}[theorem]{Remark}

\theoremstyle{break}
\newtheorem{btheorem}{Theorem}[section]
\newtheorem{algorithm}[btheorem]{Algorithm}

\title{Computing cohomology intersection numbers of GKZ hypergeometric systems}

\author{Saiei-Jaeyeong Matsubara-Heo}

\begin{document}
\date{}
\maketitle

\begin{abstract}
In this review article, we report on some recent advances on the computational aspects of cohomology intersection numbers of GKZ systems developed in \cite{GM}, \cite{MH}, \cite{MT} and \cite{MT2}. We also discuss the relation between intersection theory and evaluation of an integral of a product of powers of absolute values of polynomials.
\end{abstract}

\section{Introduction}


\subsection{Algebraic de Rham cohomology group and cohomology intersection form}
Hypergeometric functions appear in various contexts of pure and applied mathematics. Among others, Gau\ss' hypergeometric function defined by analytic continuations of ${}_2F_1\left(\substack{\alpha,\beta\\ \gamma};z\right)=\sum_{n=0}^\infty\frac{(\alpha)_n(\beta)_n}{(\gamma)_nn!}z^n$ is presumably the best studied example of a special function. Though it enjoys several properties, it is the fact that Gau\ss' hypergeometric function admits an integral representation that provides a unified means of performing analytic continuations. Reversing the perspective, one can define a hypergeometric function by means of an integral of the form
\begin{equation}\label{eqn:1stInt}
I(z)=\int_{\Gamma}\Phi\omega,
\end{equation}
where $\Phi=\prod_{l=1}^kf_l(x;z)^{\alpha_l}$, $f_l(x;z)$ is a family of polynomials in $x=(x_1,x_2,...)$ parametrized by $z=(z_1,z_2,....)$, $\alpha_l$ is a complex parameter, $\Gamma$ is a suitable integration contour, and $\omega$ is a rational top-dimensional differential form in $x$ having at most poles along $\bigcup_{l}\{ x\mid f_l(x;z)=0\}$. Integral representation (\ref{eqn:1stInt}) can be regarded as a pairing between a homology class $[\Gamma]$ and a cohomology class $[\omega]$. Therefore, the theory of algebraic de Rham cohomology groups naturally comes into play. 

We set $V_z=\{ x\mid f_1(x;z)\dots f_k(x;z)\neq 0\}$ and $\nabla_x=d_x+d_x\log\Phi\wedge$. In view of Stokes' theorem, it is natural to regard $[\omega]$ as an element of the top-dimensional de Rham cohomology group $\Homo_{dR}^{top}(V_z;\nabla_x)=\Homo^{top}(\Omega^\bullet(V_z^{alg}),\nabla_x)$ where $\Omega^\bullet(V_z^{alg})$ is the set of rational differential forms having at most poles along $\bigcup_{l}\{ x\mid f_l(x;z)=0\}$. The cohomological point of view is particularly useful when we derive the Pfaffian system of the integral (\ref{eqn:1stInt}). We consider a basis $\{[\omega_i]\}_i$ of $\Homo_{dR}^{top}(V_z;\nabla_x)$ and set $Y={}^t(I_1,I_2,...)$ where $I_i$ is defined by (\ref{eqn:1stInt}) with $\omega$ replaced by $\omega_i$. For a matrix $A(z)=(a_{ij}(z))_{i,j}$ with entries in rational 1-forms in $z$, we call the system of linear partial differential equations $dY(z)=A(z)Y(z)$ a Pfaffian system. The computation of the matrix $A(z)$ can be formalized using the Gau\ss-Manin connection $\nabla^{GM}:=d_z+d_z\log\Phi\wedge$. Taking the exterior derivative of (\ref{eqn:1stInt}), we obtain an identity $dI(z)=\int_{\Gamma}\Phi\nabla^{GM}\omega$. Therefore, the entries $a_{ij}$ are characterized by the relation $[\nabla^{GM}\omega_i]=\sum_{j}a_{ij}[\omega_j]$ in the top-dimensional de Rham cohomology group.

Another advantage of introducing cohomological point of view is that we can  relate cohomological invariants to the analysis of the integral (\ref{eqn:1stInt}). In this paper, we will focus on a particular class of invariants: {\it the cohomology intersection form}. This is a perfect bilinear pairing
\begin{equation}\label{eqn:CohIF}
\langle\bullet,\bullet\rangle_{ch}:\Homo_{dR}^{top}(V_z;\nabla_x)\times\Homo_{dR}^{top}(V_z;\nabla_x^\vee)\rightarrow\C.
\end{equation}
Here, $\nabla_x^\vee$ is the dual connection of $\nabla_x$. The value of the cohomology intersection form for a  given pair of cohomology classes is called {\it the cohomology intersection number}. Note that, in our context, we assume that the parameters $\alpha_i$ are generic.\footnote{When the parameters take special values, we should replace the de Rham cohomology groups by the {\it middle} cohomology groups. This aspect is not discussed in this paper.} The importance of this invariant in the context of hypergeometric functions was discovered in \cite{CM} by Koji Cho and Keiji Matsumoto. They showed that the cohomology intersection number naturally appears in a class of functional identities called {\it Riemann-Hodge bilinear relation.} They also developed a method of evaluating cohomology intersection numbers for algebraic de Rham cohomology groups associated to $1$-dimensional integrals. We call their method {\it residue method} because it is based on residue calculus. Residue method was later generalized to generic hyperplane arrangement case in \cite{MatsumotoIntersection}. 

Despite the fact that there have been numerous efforts to evaluate the cohomology intersection numbers (\cite{Goto}, \cite{Goto3}, \cite{Goto-Matsumoto}, \cite{MM}, \cite{MimachiYoshida},  \cite{MizeraLong}, \cite{Ohara-Sugiki-Takayama}, \cite{Weinzierl} and references therein), most of the existing methods utilize residue method in the spirit of \cite{MatsumotoIntersection}. Note, however, that in recent papers \cite{MizeraShort} and \cite{MP} the authors established a new method of evaluating the cohomology intersection numbers by means of {\it higher residues around critical points}.

The aim of this paper is to give an overview of yet another method of computing cohomology intersection numbers associated to GKZ systems  developed in \cite{MT} and \cite{MT2} based on \cite{GM} and \cite{MH}. The crucial novelties are a new characterization of the cohomology intersection form, an extensive use of {\it computational algebraic analysis}, and combinatorial structure of cohomology intersection numbers.

\subsection{The secondary equation}
Let $[\omega(z)]\in\Homo^{top}_{dR}(V_z;\nabla_x)$ and $[\omega^\vee(z)]\in\Homo^{top}_{dR}(V_z;\nabla_x^\vee)$ be cohomology classes depending rationally on $z$. Then the cohomology intersection number \newline $\langle [\omega(z)], [\omega^\vee(z)]\rangle_{ch}$ is again a rational function in $z$. The definition of the cohomology intersection form immediately gives rise to the following identity
\begin{equation}\label{eqn:2nd}
d_z\langle [\omega(z)], [\omega^\vee(z)]\rangle_{ch}=\langle \nabla^{GM}[\omega(z)], [\omega^\vee(z)]\rangle_{ch}+\langle [\omega(z)], \nabla^{GM\vee}[\omega^\vee(z)]\rangle_{ch}.
\end{equation}
Here, we have set $\nabla^{GM\vee}:=d_z-d_z\log\Phi\wedge$. We call the identity (\ref{eqn:2nd}) the {\it secondary equation}.

Let us make the secondary equation more explicit. We take bases of de Rham cohomology groups $\{ [\omega_i(z)]\}_{i=1}^r\subset\Homo^{top}_{dR}(V_z;\nabla_x)$, $\{ [\omega^\vee_i(z)]\}_{i=1}^r\subset\Homo^{top}_{dR}(V_z;\nabla_x^\vee)$ depending rationally in $z$. We trivialize the Gau\ss-Manin connections as $\nabla^{GM}=d_z+\Omega\wedge$ and $\nabla^{GM\vee}=d_z+\Omega^\vee\wedge$ with respect to these bases. Then the secondary equation (\ref{eqn:2nd}) is equivalent to the following Pfaffian system for the  {\it cohomology intersection matrix} $I_{ch}=(\langle[\omega_i(z)],[\omega^\vee_j(z)]\rangle_{ch})_{i,j=1}^r$:
\begin{equation}\label{eqn:2ndary}
d_zI_{ch}={}^t\Omega I_{ch}+I_{ch}\Omega^\vee.
\end{equation}
Thus, the cohomology intersection matrix is a rational solution of the secondary equation (\ref{eqn:2ndary}). The point is that a partial convers is also true: any rational solution of (\ref{eqn:2ndary}) is equal to the cohomology intersection matrix $I_{ch}$ up to a constant multiplication. Therefore, we can essentially evaluate any cohomology intersection number by finding a non-zero rational solution of (\ref{eqn:2ndary}).

\subsection{The viewpoint of GKZ system}

From this subsection, we will focus on a more specific integral
\begin{equation}\label{eqn:EInt}
I(z)=\int_\Gamma\prod_{l=1}^kh_l(x;z)^{-\gamma_l}x^c\omega
\end{equation}
where $h_l(x;z)=\sum_jz_j^{(l)}x^{{\bf a}^{(l)}(j)}$ are Laurent polynomials, $x^c=x_1^{c_1}x_2^{c_2}\cdots$ and $\gamma_l$ and $c_j$ are parameters. The integral (\ref{eqn:EInt}) is naturally a solution of a {\it GKZ system}, a class of holonomic systems introduced by I.M.Gelfand, M.I.Graev, M.M.Kapranov, and A.V.Zelevinsky (\cite{GGZ}, \cite{GKZToral}). Therefore, it is natural to expect that we can study the algebraic de Rham cohomology group associated to the integral (\ref{eqn:EInt}) by means of GKZ system.

The deformation parameters $z=(z_j^{(l)})_{j,l}$ can be regarded as a variable of a torus $(\C^*)^N$. We write $\DD_{\C^N}$ for the ring of linear partial differential operators on $\C^N$ with polynomial coefficients. The GKZ system is defined as a quotient $\DD$-module $M_{GKZ}:=\DD_{\C^N}/I_{GKZ}$ where $I_{GKZ}$ is a certain left ideal of $\DD_{\C^N}$. Since $M_{GKZ}$ is always holonomic (\cite{Adolphson}), it defines an integrable connection on a Zariski open dense subset $U$ of $\C^N$. The importance of GKZ system in our context is that $M_{GKZ}$ is canonically isomorphic to the Gau\ss-Manin connection on $U$. Let us consider a sheaf $\HdR^{top}$ on $U$ whose stalk $\mathcal{H}_{dR,z}^{top}$ at each $z\in U$ is canonically isomorphic to $\Homo^{top}_{dR}(V_z;\nabla_x)$ (for the precise definition, see \cite{MT}). The Gau\ss-Manin connection $\nabla^{GM}$ naturally acts on this sheaf $\HdR^{top}$ and the pair $(\HdR^{top},\nabla^{GM})$ is an integrable connection on $U$ which is canonically isomorphic to (the restriction of) $M_{GKZ}$ (\cite{GKZEuler}).

Through the isomorphism $(\HdR^{top},\nabla^{GM})\simeq M_{GKZ}$, any cohomology class $[\omega]\in\HdR^{top}$ corresponds to a modulo class $[P]\in M_{GKZ}$ represented by an operator $P\in\DD_{\C^N}$. Therefore, computations in the algebraic de Rham cohomology group $\Homo^{top}_{dR}(V_z;\nabla_x)$ are reduced to those in $M_{GKZ}$ where we can employ a toolkit of computational algebraic analysis. For example, we can compute a basis of $\Homo^{top}_{dR}(V_z;\nabla_x)$ at a generic point by computing a $\C(z)$-basis of $\C(z)\otimes_{\C[z]}M_{GKZ}$ which is equal to the set of standard monomials with respect to a Gr\"obner basis of GKZ ideal $I_{GKZ}$ for a monomial order (\cite[6.2]{dojo-en}). Once a basis of the algebraic de Rham cohomology group $\{ [\omega_1],[\omega_2],\dots\}$ is given, the connection matrix is obtained by a ``division'' of $\nabla^{GM}[\omega_i]$ by cohomology classes $[\omega_j]$. A refined version of this argument is illustrated in \S\ref{subsec:5.1}.

GKZ system also enjoys a special combinatorics, from which we can derive a formula of cohomology intersection numbers. The definition domain $(\C^*)^N$ of GKZ system admits a natural (relative) toric compactification $X$(\cite[Chapter 7]{GKZbook}, \cite{DeLoeraRambauSantos}). A remarkable fact is that at each torus fixed point of $X$, the cohomology intersection number is expanded into a convergent Laurent series whose coefficients are determined combinatorially (\cite[Theorem 8.1]{GM} and \cite[Theorem 2.6]{MH}). This is summarized in \S\ref{sec:4}.

We are in a position to illustrate how our algorithm works:

\begin{algorithm}[A prototype of the main algorithm]\ \\ 
Input: bases $\{ [\omega_i(z)]\}_{i=1}^r\subset\Homo^{top}_{dR}(V_z;\nabla_x)$, $\{ [\omega^\vee_i(z)]\}_{i=1}^r\subset\Homo^{top}_{dR}(V_z;\nabla_x^\vee)$ rational in $z$. \\
Output: the cohomology intersection matrix $I_{ch}=(\langle[\omega_i(z)],[\omega^\vee_j(z)]\rangle_{ch})_{i,j=1}^r$. 
\begin{enumerate}
\item Find connection matrices $\Omega$ and $\Omega^\vee$ of $\nabla^{GM}$ and $\nabla^{GM\vee}$ with respect to bases $\{ [\omega_i(z)]\}_{i=1}^r$ and $\{ [\omega_i^\vee(z)]\}_{i=1}^r$.
\item Find a non-zero matrix $I$ whose entries are rational functions on $U$ and which satisfies the secondary equation (\ref{eqn:2ndary}).
\item There is a complex number $C$ so that the equality $I_{ch}=C\cdot I$ holds. Specify $C$ by means of \cite[Theorem 8.1]{GM} or \cite[Theorem 2.6]{MH}. 
\end{enumerate}
\end{algorithm}

\noindent
As for step 2, we can utilize, e.g., the Maple package ``IntegrableConnections'' (\cite{IC-proj}) whose algorithm is based on \cite{Barkatou} (see also \cite{Oaku-Takayama-Tsai}). The algorithm is implemented in the computer algebra system Risa/Asir (\cite{risa-asir}, \cite{tc-web}). Combining this algorithm with the one of computing a basis of the de Rham cohomology group \cite{HNT}, we obtain a complete algorithm of determining the cohomology intersection form.

\subsection{An integral of a product of powers of absolute values of polynomials}

In addition to the algorithmic aspects, we will also discuss an integral of a product of powers of absolute values of polynomials in the last section. Namely, we consider an integral of the form
\begin{equation}\label{eqn:AbsoluteEuler}
I(\alpha)=\int_{\C^n}|\Phi|^{2}\omega\wedge\bar{\eta}.
\end{equation}
We regard (\ref{eqn:AbsoluteEuler}) as a function of parameters $\alpha_i$. This integral can be seen as a single-valued version of the integral (\ref{eqn:1stInt}) and has been studied by several people: the one-dimensional case of (\ref{eqn:AbsoluteEuler}) was discussed in \cite{HY} and a multidimensional case with specific choices of $f_l$ had appeared in \cite{MizeraLong}. 

It is classically known that $I(\alpha)$ is a meromorphic function in $\alpha=(\alpha_1,\dots,\alpha_k)\in\C^k$. The poles of $I(\alpha)$ is, in principle, described by the multivariate $b$-functions (\cite{Bahloul},\cite{Gyoja},\cite{Sabbah}). However, it is a difficult task to compute the multivariate $b$-functions in general, neither is it straightforward to obtain a closed form of the analytic continuations of $I(\alpha)$. We introduce the perspective of intersection theory to overcome this difficulty.

The important point to note here is that the integral (\ref{eqn:AbsoluteEuler}) is a variant of the cohomology intersection number. Indeed, M.Hanamura and M.Yoshida has already pointed out this fact (\cite{HY}) when $n=1$. They showed that the integral (\ref{eqn:AbsoluteEuler}) appears naturally as a polarization of $L^2$-cohomology groups. In this paper, we discuss the higher dimensional case. The basic ingredient is the theory of harmonic forms developed in \cite{KKPoincare}. By writing down the Riemann-Hodge bilinear relation in this context, we obtain a method of computing the analytic continuation of the integral (\ref{eqn:AbsoluteEuler}). This is achieved in \S\ref{sec:61}. When the integrand $\Phi$ is related to GKZ system, we obtain a series expansion of (\ref{eqn:AbsoluteEuler}) in terms of hypergeometric series in\S\ref{sec:62}. It is expected that the theory of $b$-functions is related to our approach from intersection theory. We will not address this problem in this paper.

Finally, we remark that we do not give proofs of the statements in this paper. The proofs are available in \cite{GM}, \cite{MH}, \cite{MT} and \cite{MT2} except for the results in \S6. A more comprehensive treatment of the results in \S6 will appear elsewhere.

\section{Basic set-ups}\label{sec:Basic}
This section is devoted to recalling basic notions and notation related to algebraic de Rham cohomology groups. A more comprehensive description can be found in \cite{AK} or in \cite{Deligne}. The readers familiar with these notions can skip this section.

\subsection{Algebraic de Rham cohomology groups}
We fix a positive integer $n$ and consider non-constant complex polynomials $f_l(x)$ $(l=1,\dots,k)$ in $x\in\C^n$. We choose complex numbers $\alpha_l\in\C$ $(l=1,\dots,k)$ and set $\Phi:=\prod_{l=1}^kf_l(x)^{\alpha_l}$. We are interested in the integral of the form
\begin{equation}\label{eqn:IntroInt}
I=\int_\Gamma\Phi\omega
\end{equation}
where $\Gamma$ is a suitable cycle and $\omega$ is an algebraic $n$-form in $x$. In order for the integral (\ref{eqn:IntroInt}) to define a function, we need to introduce a deformation variable $z$ in $f_l$, namely we consider the case when $f_l=f_l(x;z)$ depends polynomially in $z$ and therefore $I=I(z)$ is an analytic function in $z$. Under the presence of $z$, the integral $I(z)$ is called an Euler integral representation. For the moment, we fix the deformation variable $z$ to make the dependence on $z$ implicit. We set $V=\{ x\in\C^n\mid f_1(x)\cdots f_k(x)\neq 0\}$. The twisted differential $\nabla_x$ associated to the integral (\ref{eqn:IntroInt}) is defined by $d_x+\sum_{l=1}^k\alpha_ld_x\log f_l\wedge$ where $d_x$ is the exterior derivative on $V$. Let us denote by $\Omega^p(V^{alg})$ the set of algebraic differential $p$-forms on $V$. It can readily be seen that the sequence
\begin{equation}\label{eqn:complex1}
\cdots\rightarrow\Omega^p(V^{alg})\overset{\nabla_x}{\rightarrow}\Omega^{p+1}(V^{alg})\rightarrow\cdots
\end{equation}
defines a complex where $\Omega^p(V^{alg})$ is at the p-th position. We define $p$-th algebraic de Rham cohomology group $\Homo^p_{dR}\left( V^{alg};\nabla_x\right)$ as the $p$-th cohomology group of the complex (\ref{eqn:complex1}). We use the symbol $V^{an}$ to emphasize that we equip $V$ with the analytic topology in contrast to the same set $V^{alg}$ equipped with Zariski topology. Since we are mainly interested in the algebraic de Rham cohomology group rather than the analytic one (though they are isomorphic), we simply write $\Homo^p_{dR}\left( V;\nabla_x\right)$ for $\Homo^p_{dR}\left( V^{alg};\nabla_x\right)$.

We define the dual object of $\Homo^p_{dR}\left( V;\nabla_x\right)$. We write  $\mathcal{L}$ for the dual local system of flat sections of $\nabla_x^{an}$. More intuitively, any local section of the sheaf $\mathcal{L}$ is a complex number times a determination of the multivalued function $\Phi$. With this notation, we can define $p$-th {\it twisted homology group} $\Homo_p\left( V^{an};\mathcal{L}\right)$ (\cite{Steenrod}). For readers' convenience, we give an explicit description of the twisted homology group.

 We write $C_p(V^{an};\mathcal{L})$ for the vector space of formal finite sums $\Gamma=\sum_ia_i\Gamma_i\otimes \Phi|_{\Gamma_i}$ where $a_i\in\C$, $\Gamma_i$ is a continuous map $\Gamma_i:\Delta^p\rightarrow V^{an}$ and $\Phi|_{\Gamma_i}$ is a determination of the multivalued function $\Phi$ on the image of $\Gamma_i$. Here, the symbol $\Delta^p$ stands for the $p$-dimensional simplex. Let us denote by $\partial\Gamma_i$ the boundary of $\Gamma_i$ in the ordinary sense. Setting $\partial_\Phi\Gamma=\sum_ia_i\partial\Gamma_i\otimes \Phi|_{\partial\Gamma_i}$, the sequence
\begin{equation}\label{eqn:complex2}
\cdots\rightarrow C_p(V^{an};\mathcal{L})\overset{\partial_\Phi}{\rightarrow}C_{p-1}(V^{an};\mathcal{L})\rightarrow\cdots
\end{equation}
defines a complex where $C_p(V^{an};\mathcal{L})$ is at the $(-p)$-th position. We define $p$-th twisted homology group $\Homo_p\left( V^{an};\mathcal{L}\right)$ as the $(-p)$-th cohomology group of the complex (\ref{eqn:complex2}). For any twisted $p$-chain $\Gamma=\sum_ia_i\Gamma_i\otimes\Phi\in C_p(V^{an};\mathcal{L})$ and an algebraic $p$-form $\omega\in\Omega^p(V^{alg})$, we set $\langle\Gamma,\omega\rangle_{per}=\sum_ia_i\int_{\Gamma_i}\Phi\omega$. Note that the determination of $\Phi$ is specified on $\Gamma_i$. Then, it is classical that the pairing 
\begin{equation}\label{eqn:Period}
\begin{array}{cccc}
\langle\bullet,\bullet\rangle_{per}:&\Homo_p\left( V^{an};\mathcal{L}\right)\times \Homo^p_{dR}\left( V;\nabla_x\right)&\rightarrow&\mathbb{C}\\
&\rotatebox{90}{$\in$}& &\rotatebox{90}{$\in$}\\
&([\Gamma],[\omega])&\mapsto&\ \ \ \langle\Gamma,\omega\rangle_{per}
\end{array}
\end{equation}
is well-defined and gives rise to a perfect pairing (\cite{Deligne}).

We write $\mathcal{L}^\vee$ for the dual local system of $\mathcal{L}$, that is, any local section of $\mathcal{L}^\vee$ is a complex number times a determination of $\Phi^{-1}$. In the same way, we write $\nabla_x^\vee$ for the connection $d_x-\sum_{l=1}^k\alpha_ld_x\log f_l\wedge$ which is dual to $\nabla_x$. We can also define the perfect bilinear pairing 
\begin{equation}\label{eqn:Period2}
\begin{array}{cccc}
\langle\bullet,\bullet\rangle_{per}:&\Homo_p\left( V^{an};\mathcal{L}^\vee\right)\times \Homo^p_{dR}\left( V;\nabla_x^\vee\right)&\rightarrow&\mathbb{C}\\
&\rotatebox{90}{$\in$}& &\rotatebox{90}{$\in$}\\
&([\Gamma^\vee],[\omega^\vee])&\mapsto&\ \ \ \langle\Gamma^\vee,\omega^\vee\rangle_{per}
\end{array}
\end{equation}
where $\langle\Gamma^\vee,\omega^\vee\rangle_{per}$ is defined as an integration as in (\ref{eqn:Period}).

\subsection{The (co)homology intersection form}

We write $\Homo^n_{dR,c}\left( V^{an};\nabla_x^{an}\right)$ for the analytic de Rham cohomology group with compact support. Namely, if the symbol $\mathcal{E}^p_c(V^{an})$ denotes the set of smooth $p$-forms on $V^{an}$ with compact support, $\Homo^n_{dR,c}\left( V^{an};\nabla_x^{an}\right)$ is defined as the $n$-th cohomology group of the complex $(\mathcal{E}^\bullet_c(V^{an}),\nabla_x^{an})$. By Poincar\'e-Verdier duality, the bilinear pairing 
\begin{equation}
\begin{array}{ccc}
\Homo^n_{dR,c}\left( V^{an};\nabla_x^{an}\right)\times\Homo^n_{dR}\left( V^{an};\nabla_x^{an\vee}\right)&\rightarrow&\mathbb{C}\\
\rotatebox{90}{$\in$}& &\rotatebox{90}{$\in$}\\
([\omega],[\omega^\vee])&\mapsto&\int_{V^{an}}\omega\wedge\omega^\vee
\end{array}
\end{equation}
is perfect. Let $\mathcal{E}^p(V^{an})$ be the set of smooth $p$-forms on $V^{an}$. The natural inclusion $\mathcal{E}^p_c(V^{an})\hookrightarrow\mathcal{E}^p(V^{an})$ induces a morphism of cohomology groups $can:\Homo^p_{dR,c}\left( V^{an};\nabla_x^{an}\right)$ $\rightarrow$ $\Homo^p_{dR}\left( V^{an};\nabla_x^{an}\right)$. We say that {\it the regularization condition} is satisfied if the morphism $can$ is an isomorphism for any $p$. Note that the regularization condition is a generic condition for the parameters $\alpha_l$. The regularization condition implies the pure-codimensionality of the cohomology groups. Namely, we have the vanishing $\Homo^p_{dR}\left( V^{an};\nabla_x^{an}\right)=0$ for any $p\neq n$. In the following, we always assume that the regularization condition is satisfied.\footnote{It is also important to study the integral (\ref{eqn:IntroInt}) when the regularization condition is violated. In this case, the subsequent argument can be developed in a parallel way if we replace the algebraic de Rham cohomology group by the so-called {\it middle} cohomology group and replace each perfect pairing by its middle version (\cite{FSY}).  However, it seems that there is no systematic way of computing middle cohomology groups by means of computer algebra. This is the future task.} A criterion for this assumption is explained in \S\ref{GKZsystem}. Since $\nabla_x$ is a regular connection, the canonical morphism $\Homo^n_{dR}\left( V;\nabla_x\right)\rightarrow\Homo^n_{dR}\left( V^{an};\nabla_x^{an}\right)$ is always an isomorphism by Deligne-Grothendieck comparison theorem (\cite[Corollaire 6.3]{Deligne}). Therefore, we have a canonical isomorphism ${\rm reg}:\Homo^n_{dR}\left( V;\nabla_x\right)\rightarrow\Homo^n_{dR,c}\left( V^{an};\nabla_x^{an}\right)$. Note that the Poincar\'e dual of the isomorphism ${\rm reg}$ is called a regularization map in the theory of special functions (\cite[\S 3.2]{Aomoto-Kita}). Finally, we define the cohomology intersection form $\langle\bullet,\bullet\rangle_{ch}$ between algebraic de Rham cohomology groups by the formula
\begin{equation}\label{eqn:19}
\begin{array}{cccc}
\langle\bullet,\bullet\rangle_{ch}:&\Homo^n_{dR}\left( V;\nabla_x\right)\times\Homo^n_{dR}\left( V;\nabla_x^{\vee}\right)&\rightarrow&\mathbb{C}\\
&\rotatebox{90}{$\in$}& &\rotatebox{90}{$\in$}\\
&([\omega],[\omega^\vee])&\mapsto&\int_{V^{an}}{\rm reg}([\omega])\wedge\omega^\vee.
\end{array}
\end{equation}
The value $\langle[\omega],[\omega^\vee]\rangle_{ch}$ is called {\it the cohomology intersection number} of
$[\omega]$ and $[\omega^\vee]$.

\begin{rem}
In the definition (\ref{eqn:19}) of the cohomology intersection form, we can also obtain the same perfect pairing by regularizing the form $\omega^\vee$. Namely, we have equalities
\begin{equation}
\langle[\omega],[\omega^\vee]\rangle_{ch}=\int_{V^{an}}{\rm reg}([\omega])\wedge\omega^\vee=\int_{V^{an}}\omega\wedge{\rm reg}([\omega^\vee])=\int_{V^{an}}{\rm reg}([\omega])\wedge{\rm reg}([\omega^\vee]).
\end{equation}
\end{rem}

In view of the perfect pairings (\ref{eqn:Period}) and (\ref{eqn:19}), we can define a natural perfect bilinear pairing between twisted homology groups

\begin{equation}\label{eqn:HIF}
\langle\bullet,\bullet\rangle_{h}:\Homo_n\left( V^{an};\mathcal{L}\right)\times \Homo_n\left( V^{an};\mathcal{L}^\vee\right)\rightarrow\mathbb{C}
\end{equation}
which we call the {\it homology intersection form}. It is known that one can compute the homology intersection number by counting the geometric intersections. We take twisted cycles $[\Gamma]\in\Homo_n\left( V^{an};\mathcal{L}\right)$ and $[\Gamma^\vee]\in\Homo_n\left( V^{an};\mathcal{L}^\vee\right)$. If $[\Gamma]$ (resp. $[\Gamma^\vee]$) is represented by a chain $\sum_ia_i\Gamma_i\otimes\Phi$ (resp. $\sum_ia_i^\prime\Gamma_i^\prime\otimes\Phi^{-1}$), the intersection number $\langle [\Gamma],[\Gamma^\vee]\rangle_h$ is equal to $\sum_{i,j}a_ia_j^\prime I_{loc}(\Gamma_i,\Gamma_j^\prime)$. Here, $I_{loc}(\Gamma_i,\Gamma_j^\prime)$ is the local intersection number (see e.g. \cite{KY}).

\subsection{Twisted period relations}

It was discovered in \cite{CM} that a family of functional identities of hypergeometric functions called quadratic relations can be derived in a systematic way from the Riemann-Hodge bilinear relation. This relation is a compatibility among cohomology intersection form (\ref{eqn:19}), homology intersection form (\ref{eqn:HIF}), and period pairings (\ref{eqn:Period}) and  (\ref{eqn:Period2}). Let us take four bases $\{[\omega_i]\}_{i=1}^r\subset\Homo_{dR}^n(V;\nabla_x)$, $\{[\omega_i^\vee]\}_{i=1}^r\subset\Homo_{dR}^n(V;\nabla_x^\vee)$, $\{[\Gamma_i]\}_{i=1}^r\subset\Homo_n(V^{an};\mathcal{L})$, $\{[\Gamma_i^\vee]\}_{i=1}^r\subset\Homo_n(V^{an};\mathcal{L}^\vee)$. We set $I_{ch}:=\left(\langle [\omega_i],[\omega_j^\vee]\rangle_{ch}\right)_{i,j=1}^r$, $I_{h}:=\left(\langle [\Gamma_i],[\Gamma_j^\vee]\rangle_{h}\right)_{i,j=1}^r$, $P:=\left(\langle [\Gamma_j],[\omega_i]\rangle_{per}\right)_{i,j=1}^r$, $P^\vee:=\left(\langle [\Gamma_j^\vee],[\omega_i^\vee]\rangle_{per}\right)_{i,j=1}^r$. 

\begin{theorem}[Riemann-Hodge bilinear relation, or, twisted period relation]
\begin{equation}
I_{ch}=P{}^tI_h^{-1}{}^tP^\vee.
\end{equation}
\end{theorem}
In particular, if we write $I_h^{-1}=(C^{ij})_{i,j=1}^r$, we obtain an expansion formula of the cohomology intersection number
\begin{equation}\label{eqn:expansion}
\langle [\omega],[\omega^\vee]\rangle_{ch}=\sum_{i,j=1}^r\left(\int_{\Gamma_i}\Phi\omega\right)C^{ji}\left(\int_{\Gamma_j^\vee}\Phi^{-1}\omega^\vee\right)
\end{equation}
for any cohomology classes $[\omega]\in\Homo_{dR}^n(V;\nabla_x)$ and $[\omega^\vee]\in\Homo_{dR}^n(V;\nabla_x^\vee)$. In view of the formula (\ref{eqn:expansion}), we can evaluate the cohomology intersection number in terms of the periods and homology intersection numbers. We will see that these functions can explicitly be evaluated when the integral comes from a GKZ system in \S\ref{sec:4}.

\section{Gau\ss-Manin connection and the secondary equation} \label{sec:2}

\subsection{Gau\ss-Manin connection}
Recall that the polynomials $f_l=f_l(x;z)$ in \S\ref{sec:Basic} depend implicitly on other complex variables $z=(z_1,\dots,z_N)$. Therefore, we suppose there are smooth affine varieties $U,V$ and a smooth affine morphism $\pi:V\rightarrow U$ so that each fiber is given by $V_z:=\pi^{-1}(z)=\{ x\in\C^n\mid f_1(x;z)\cdots f_k(x;z)\neq 0\}$. Note that $f_l$ can be regarded as regular functions on $V$. For each fiber $V_z$, we can carry out the construction in the previous section to define the algebraic de Rham cohomology groups, twisted homology groups, and perfect pairings among them. It is natural to expect that the algebraic de Rham cohomology group $\Homo^n_{dR}(V_z;\nabla_x)$ depends rationally on the variables $z$. This is the viewpoint of Gau\ss-Manin connection.   

Let us formulate the Gau\ss-Manin connection briefly. As for the rigorous description in the present context, see \cite{MT}. We regard $f_l$ as a global section of the sheaf of regular functions $\mathcal{O}_V$ on $V$ and set $\nabla=d+\sum_{l=1}^k\alpha_ld\log f_l\wedge$. Here, $d$ is the exterior derivative on $V$. We assume that the morphism $\pi$ is locally given by a projection and the dimension of each fiber $V_z$ is $n$. Therefore, we can decompose $V$ as a product of the base space $U$ and the fiber $V_z$ for some $z$. Let $x$ be a coordinate of the fiber $V_z$. We can write $\nabla$ as a sum $\nabla=\nabla_x+\nabla_z$, where $\nabla_x$ (resp. $\nabla_z$) is defined by $d_x+\sum_{l=1}^k\alpha_ld_x\log f_l\wedge$ (resp. $d_z+\sum_{l=1}^k\alpha_ld_z\log f_l\wedge$). We define the relative de Rham cohomology group $\mathcal{H}^n_{dR}$ as the $n$-th cohomology group of the complex
\begin{equation}
\left(0\rightarrow\Omega^0_{V/U}\overset{\nabla_x}{\rightarrow}\Omega^1_{V/U}\overset{\nabla_x}{\rightarrow}\cdots\overset{\nabla_x}{\rightarrow}\Omega^n_{V/U}\rightarrow 0\right).
\end{equation}
Here, $\Omega^p_{V/U}$ denotes the sheaf of relative differential $p$-forms $\omega(z)$ locally defined by $\sum_{|I|=p}a(x;z)dx^{I}$ ($a(x;z)\in \mathcal{O}_V$). For any $z\in U$, there is a natural evaluation morphism ${\rm ev}_z:\mathcal{H}_{dR,z}^n\rightarrow\Homo^n_{dR}\left( V_z;\nabla_x\right)$. We define the dual object $\HdR^{n\vee}$ by replacing $\nabla_x$ by $\nabla_x^\vee$ in the construction above.  Therefore, for any local section $[\omega]$ of $\HdR^n$ and $[\omega^\vee]$ of $\HdR^{n\vee}$, we can define the cohomology intersection number $\langle[\omega],[\omega^\vee]\rangle_{ch}$ as a function of $z\in U$ by the formula $U\ni z\mapsto\langle{\rm ev}_z([\omega]),{\rm ev}_z([\omega^\vee])\rangle_{ch}\in\C$. This actually defines a $\mathcal{O}_U$-bilinear morphism $\langle\bullet,\bullet\rangle_{ch}:\HdR^n\times\HdR^{n\vee}\rightarrow\mathcal{O}_U$.

We define the Gau\ss-Manin connection $\nabla^{GM}:\HdR^n\rightarrow\Omega^1_U(\HdR^n):=\Omega_U^1 \otimes \HdR^n$. For any local section $[\omega]$ of $\HdR^n$, we set
\begin{equation}\label{GaussManin}
\nabla^{GM}[\omega]:=[\nabla_z\omega].
\end{equation}
Here, the superscript $GM$ stands for "Gau\ss-Manin". The operator $\nabla^{GM}$ may change the cohomology class but the result does not depend on a particular choice of representatives. The dual connection $\nabla^{GM\vee}:\HdR^{n\vee}\rightarrow\Omega^1_U(\HdR^{n\vee})$ is defined by replacing $\nabla_z$ by $\nabla_z^\vee$ in (\ref{GaussManin}).

\subsection{The secondary equation}

We can rewrite the action of $\nabla^{GM}$ in terms of local frames. Let $\{ [\omega_i]\}_{i=1}^r$ be a local free basis of $\HdR^n$. We define the connection matrix $\Omega=(\omega_{ij})_{i,j=1}^r$ whose entries are $1$-forms on $U$ so that we have an equality $\nabla^{GM}[\omega_i]=\sum_{j=1}^r\omega_{ji}\wedge[\omega_j]$. Then, the Gau\ss-Manin connection $\nabla^{GM}$ is given by $\nabla^{GM}=d_z+\Omega\wedge$ with respect to the local frame $\{ [\omega_i]\}_{i=1}^r$. In the same way, we can also define the connection matrix $\Omega^\vee$ of $\nabla^{GM\vee}$ with respect to the local frame $\{ [\omega_i^\vee]\}_{i=1}^r$ of $\HdR^{n\vee}$. These connection matrices can be seen as Pfaffian systems. We take cycles $[\Gamma]\in\Homo_n(X_z^{an};\mathcal{L})$ and $[\Gamma]\in\Homo_n(X_z^{an};\mathcal{L}^\vee)$). We put $Y={}^t\left(\int_{\Gamma}\Phi\omega_1,\dots,\int_{\Gamma}\Phi\omega_r\right)$ and $Y^\vee={}^t\left(\int_{\Gamma^\vee}\Phi^{-1}\omega^\vee_1,\dots,\int_{\Gamma^\vee}\Phi^{-1}\omega^\vee_r\right)$. Then, we have Pfaffian systems 

\begin{equation}
d_z Y={}^t\Omega Y \ \ \text{      and      }\ \  d_z Y^\vee={}^t\Omega^\vee Y^\vee.
\end{equation}

The Gau\ss-Manin connections $\nabla^{GM}$ and $\nabla^{GM\vee}$ on $\HdR^n$ and $\HdR^{n\vee}$ are compatible with the cohomology intersection form. Namely, for any local sections $[\omega]$ of $\HdR^n$ and $[\omega^\vee]$ of $\HdR^{n\vee}$, we have
\begin{equation}\label{SecondaryEq}
d_z\langle [\omega],[\omega^\vee]\rangle_{ch}=\langle \nabla^{GM}[\omega],[\omega^\vee]\rangle_{ch}+\langle [\omega],\nabla^{GM\vee}[\omega^\vee]\rangle_{ch}.
\end{equation}
We call (\ref{SecondaryEq}) the secondary equation. Let us rewrite it in terms of local frames. We set $I=I_{ch}=(\langle [\omega_i],[\omega^\vee_j]\rangle_{ch})_{i,j}$ and call it the cohomology intersection matrix. Then, the secondary equation (\ref{SecondaryEq}) is equivalent to the system 
\begin{equation}\label{SecondaryEq2}
d_zI={}^t\Omega I+I\Omega^\vee.
\end{equation}
We also call (\ref{SecondaryEq2}) the secondary equation. The theorem which our algorithm is based on is the following

\begin{theorem}\label{main} {\rm \cite{MT}}
Under the regularization condition, all the entries of the cohomology intersection matrix $I_{ch}$ are rational functions. Moreover, any rational function solution $I$ of the secondary equation (\ref{SecondaryEq2}) is, up to a scalar multiplication, equal to $I_{ch}$.
\end{theorem}

\section{Intersection theory and GKZ system}\label{sec:GKZ}

\subsection{GKZ systems and Euler integral representations}\label{GKZsystem}
From this section, we slightly change the notation and consider an integral of the form
\begin{equation}\label{eqn:EulerInt}
I=\int_\Gamma\prod_{l=1}^kh_l(x)^{-\gamma_l}x^c\omega
\end{equation}
where $\Gamma$ is a suitable cycle, $\omega$ is an algebraic $n$-form in $x=(x_1,\dots, x_n)$, $\gamma_l,c_i\in \C$ are parameters ($l=1,\dots,k$ $i=1,\dots,n$),  $x^c=x_1^{c_1}\cdots x_n^{c_n}$, and $h_l(x)=h_l(x;z)=\sum_{j=1}^{N_l}z_j^{(l)}x^{{\bf a}^{(l)}(j)}$ are Laurent polynomials in $x$. Hereafter, we set
\begin{equation}
\Phi=\prod_{l=1}^kh_l(x)^{-\gamma_l}x^c.
\end{equation}
The setting of \S\ref{sec:2} is now simplified as follows: we set $N=N_1+\cdots+N_k$ and write $z$ for $(z_j^{(l)})_{j,l}$. Let $V$ be the smooth affine algebraic variety defined by $V=\{ (z,x)\in\C^N\times(\C^*)^n\mid \prod_{l=1}^kh_l(x;z)\neq 0\}$ and let $\pi:V\rightarrow \C^N$ be the projection. By \cite[2.9]{GKZEuler}, we can find a Zariski open subset $U$ of $\C^N$ on which the projection $\pi$ satisfies the assumption of \S\ref{sec:2}. Therefore, we may replace $V$ by $\pi^{-1}(U)$ to assume that $\HdR^n$ is an algebraic vector bundle on $U$.   

In this setting, the Gau\ss-Manin connection $\HdR^n$ can be extended to a holonomic $\DD$-module on $\C^N$ called the GKZ system. Let us recall the definition of GKZ system (\cite{GKZToral}). For a given $d\times n$ ($d<n$) integer matrix $A=({\bf a}(1)|\cdots|{\bf a}(n))$ and a parameter vector $\delta\in\C^{d}$, GKZ system $M_A(\delta)$ is defined as a system of partial differential equations on $\C^n$ given by
\begin{subnumcases}{M_A(\delta):}
E_i\cdot f(z)=0 &($i=1,\dots, d$)\label{EulerEq}\\
\Box_u\cdot f(z)\hspace{-0.8mm}=0& $\left(u\in \Ker(A\times:\Z^{n\times 1}\rightarrow\Z^{d\times 1})\right)$,\label{ultrahyperbolic}
\end{subnumcases}
where $E_i$ and $\Box_u$ for $u={}^t(u_1,\dots,u_{n})$ are differential operators defined by 

\begin{equation}\label{HGOperators}
E_i=\sum_{j=1}^{n}a_{ij}z_j\frac{\partial}{\partial z_j}+\delta_i,\;\;\;
\Box_u=\prod_{u_j>0}\left(\frac{\partial}{\partial z_j}\right)^{u_j}-\prod_{u_j<0}\left(\frac{\partial}{\partial z_j}\right)^{-u_j}.
\end{equation}

\noindent
We write $\DD_{\C^n}$ for the ring of linear partial differential operators on $\C^n$ with polynomial coefficients. The GKZ ideal is a left ideal of $\DD_{\C^n}$ defined by $I_{GKZ}:=\DD_{\C^n}\langle E_i,\Box_u\mid i=1,\dots,d,\ u\in L_A\rangle$. As a $\DD$-module, we set $M_A(\delta):=\DD_{\C^n}/I_{GKZ}$. It is known that GKZ system $M_A(\delta)$ is holonomic (\cite{Adolphson}). For convenience, we assume an additional condition $\Z A:=\Z{\bf a}(1)+\dots+\Z{\bf a}(n)=\Z^{d}$. In our setting, we put $A_l=({\bf a}^{(l)}(1)|\dots|{\bf a}^{(l)}(N_l))$, $d=n+k$, $n=N$. We  define an $(n+k)\times N$ matrix $A$ by


\begin{equation}\label{CayleyConfigu}
A
=
\left(
\begin{array}{ccc|ccc|c|ccc}
1&\cdots&1&0&\cdots&0&\cdots&0&\cdots&0\\
\hline
0&\cdots&0&1&\cdots&1&\cdots&0&\cdots&0\\
\hline
 &\vdots& & &\vdots& &\ddots& &\vdots& \\
\hline
0&\cdots&0&0&\cdots&0&\cdots&1&\cdots&1\\
\hline
 &A_1& & &A_2& &\cdots & &A_k& 
\end{array}
\right).
\end{equation}

\noindent
We put $\delta=(\gamma_1,\dots,\gamma_k, c_1,\dots,c_n)$. We say that the parameter $\delta$ is non-resonant if it does not belong to $\C\Gamma+\Z^d$ for any facet $\Gamma$ of the cone $\sum_{j=1}^N\R_{\geq 0}{\bf a}(j)$. If the parameter vector $\delta$ is non-resonant and $\gamma_l\notin\Z$, the Gau\ss-Manin connection $(\HdR^n,\nabla^{GM})$ (resp. $(\HdR^{n\vee},\nabla^{GM\vee})$) is isomorphic to a restriction of the GKZ system $M_A(\delta)\restriction_U$ (resp. $M_A(-\delta)\restriction_U$) and the regularization condition is true (\cite[Theorem 2.12]{MH}\footnote{This was originally proved in \cite[2.9]{GKZEuler} where the condition $\gamma_l\notin\Z$ is missing. If some of $\gamma_l$ takes an integral value, the regularization condition is, in general, violated.}). We set $\frac{dx}{x}=\frac{dx_1}{x_1}\wedge\dots\wedge\frac{dx_n}{x_n}$. The isomorphism $M_A(\delta)\restriction_U\simeq\HdR^n$ is given by the correspondence $[1]\mapsto [\frac{dx}{x}]$. Thus, any section $[\omega]$ of $\HdR^n$ can be written as $[\omega]=P\cdot [\frac{dx}{x}]$ for some linear differential operator $P\in\DD_{\C^N}$. The action of a partial derivative $\partial_{j,l}:=\frac{\partial}{\partial z_j^{(l)}}$ onto a cohomology class $[\omega(z)]\in\HdR^n$ is concretely given by
\begin{equation}\label{eqn:Action}
\partial_{j,l}\cdot [\omega(z)]=\nabla^{GM}_{\partial_{j,l}}[\omega(z)]=\left( \partial_{j,l}-\gamma_l\frac{x^{{\bf a}^{(l)}(j)}}{h_l(x;z)}\right)[\omega(z)].
\end{equation}
In the last expression of (\ref{eqn:Action}), $\partial_{j,l}$ acts as a standard partial derivative to each coefficient of the differential form $\omega(z)=\sum_ia_i(x;z)dx$.

Based on this description of the isomorphism $M_A(\delta)\restriction_U\simeq\HdR^n$, we can obtain an algorithm of computing a (free) basis of the algebraic de Rham cohomology group (\cite{HNT}). 

\begin{theorem}[\cite{HNT}, Theorem 2]\label{thm:thm2}
Suppose $\delta$ is non-resonant and $\gamma_l\notin\Z$\footnote{This condition is based on the description of the isomorphism $M_A(\delta)\restriction_U\simeq\HdR^n$ (\cite{MH}).}. If $\{ [\partial^\alpha]\}_{\alpha}\subset \C(z)\otimes_{\C[z]}M_A(\delta)$ is a $\C(z)$-basis of $\C(z)\otimes_{\C[z]}M_A(\delta)$, then $\{ \partial^\alpha\cdot [\frac{dx}{x}]\}_{\alpha}$ is a basis of $\Homo_{dR}^n(V_z;\nabla_x)$ for generic $z$.
\end{theorem}
An important point of the Theorem above is that a $\C(z)$-basis of $\C(z)\otimes_{\C[z]}M_A(\delta)$ can be computed from a Gro\"obner basis with respect to a monomial order. Let $G$ be the Gr\"obner basis of the GKZ ideal $I_{GKZ}$. The set of the standard monomials for $G$,
which is the set of monomials in $\pd{}$ that are not divisible
by the elements of $G$, is of the form
$\{ \pd{}^\alpha \,|\, \alpha \in S \subset \Z_{\geq0}^N\}$  
(see, e.g., \cite[6.1]{dojo-en}). Theorem \ref{thm:thm2} implies that $\{ \partial^\alpha\cdot [\frac{dx}{x}]\}_{\alpha\in S}$ is a basis of $\Homo_{dR}^n(V_z;\nabla_x)$ for generic $z$. Moreover, the Gr\"obner basis technique also provides an algorithm of computing the Pfaffian systems, or equivalently, the connection matrix (\cite[6.2]{dojo-en}). In \S\ref{subsec:5.1}, we discuss an improved version of this algorithm.

\subsection{Combinatorics of homology intersection numbers}\label{sec:4}


From this subsection, we assume that $\delta$ is non-resonant and $\gamma_l\notin\Z$. Let $z\in\C^N$ be a generic point. To be more precise, $z$ is taken to be a nonsingular point in the sense of \cite[Definition 3.7]{MH}. The isomorphism of $\DD$-modules $M_A(\delta)\restriction_U\simeq\HdR^n$ gives rise to an isomorphism of the solution space of GKZ system and a  twisted homology group. Namely, we have a canonical isomorphism 
\begin{equation}\label{eqn:Integration}
\Homo_n(V_z^{an},\mathcal{L})\tilde{\rightarrow}\sol_{M_A(\delta),z},
\end{equation}
the correspondence of which is given by $\Homo_n(V_z^{an},\mathcal{L})\ni[\Gamma]\mapsto \int_\Gamma \Phi\frac{dx}{x}$ (\cite[Theorem 2.10]{GKZEuler}, \cite[Theorem 2.12]{MH}). Here, $\sol_{M_A(\delta),z}$ is the stalk of the solution sheaf of $M_A(\delta)$ at the point $z$. The solution space $\sol_{M_A(\delta),z}$ has a combinatorial structure when $z$ is close to a special point in a suitable toric compactification. Through the isomorphism (\ref{eqn:Integration}), we can introduce a combinatorial structure in the twisted homology group $\Homo_n(V_z^{an},\mathcal{L})$. In order to describe this combinatorial structure, we briefly recall basic definitions related to regular triangulations.


Let us recall the definition of a regular triangulation. In general, for any subset $\sigma$ of $\{1,\dots,N\},$ $\cone(\sigma)$ denotes the positive span of the column vectors of $A$ $\{{\bf a}(1),\dots,{\bf a}(N)\}$ i.e., $\cone(\sigma)=\displaystyle\sum_{i\in\sigma}\R_{\geq 0}{\bf a}(i).$ We often identify a subset $\sigma\subset\{1,\dots,N\}$ with the corresponding set of vectors $\{{\bf a}(i)\}_{i\in\sigma}$ or with the set $\cone(\s)$. A collection $T$ of subsets of $\{1,\dots,N\}$ is called a triangulation if $\{\cone(\sigma)\mid \sigma\in T\}$ is the set of cones in a simplicial fan whose support equals $\cone(A)$. We regard $\Z^{1\times N}$ as the dual lattice of $\Z^{N\times 1}$ via the standard dot product. Let $\pi_A:\Z^{1\times N}\rightarrow L_A^\vee$ be the dual morphism of the natural inclusion $L_A\hookrightarrow \Z^{N\times 1}$ where $L_A^\vee$ is the dual lattice $\Hom_{\Z}(L_A,\Z)$. By abuse of notation, we continue to write $\pi_A$ for the linear map $\pi_A\underset{\Z}{\otimes}{\rm id}_{\R}:\R^{1\times N}\rightarrow L_A^\vee\underset{\Z}{\otimes}\R$ where ${\rm id}_{\R}:\R\rightarrow\R$ is the identity map. Then, for any generic choice of a vector $\omega\in\R^{1\times N},$ we can define a triangulation $T(\omega)$ as follows: A subset $\sigma\subset\{1,\dots,N\}$ belongs to $T(\omega)$ if there exists a vector ${\bf n}\in\R^{1\times (n+k)}$ such that ${\bf n}\cdot{\bf a}(i)=\omega_i$ if $i\in\sigma$ and ${\bf n}\cdot{\bf a}(j)<\omega_j$ if $ j\notin\sigma.$ A triangulation $T$ is called a regular triangulation if $T=T(\omega)$ for some $\omega\in\R^{1\times N}.$ For a regular triangulation $T$, we set
\begin{equation}
C_T=\Big\{ \omega\in\R^{1\times N}\mid T(\omega)=T\Big\}.
\end{equation}

\noindent
We cite a fundamental result of Gelfand, Kapranov, and Zelevinsky (\cite[Theorem 5.2.11.]{DeLoeraRambauSantos},\cite[Chapter 7, Proposition 1.5.]{GKZbook}).

\begin{theorem}[\cite{DeLoeraRambauSantos},\cite{GKZbook}]
There exists a complete fan ${\rm Fan}(A)$ in $\R^{1\times N}$ whose maximal cones are precisely $\{ C_T\}_{T: \text{regular triangulation}}$. The fan ${\rm Fan}(A)$ is called the secondary fan.
\end{theorem}

\begin{rem}
Let $F$ be a fan obtained by applying the projection $\pi_A$ to each cone of ${\rm Fan}(A)$. By definition, each cone of ${\rm Fan}(A)$ is a pull-back of a cone of $F$ through the projection $\pi_A$. Therefore, the fan $F$ is also called the secondary fan.
\end{rem}


For any subset $\tau\subset\{1,\dots,N\}$, $A_\tau$ denotes the matrix given by the columns of $A$ indexed by $\tau.$ We say that a regular triangulation $T$ is unimodular if we have $\det A_\s=\pm 1$ for any simplex $\s\in T$. Though a unimodular regular triangulation may not exist in general, many interesting classes of GKZ system admit it. For example, GKZ system coming from Appell-Lauricella system or Horn's system has a unimodular regular triangulation. 

In order to simplify the exposition, we assume the matrix $A$ defined by (\ref{CayleyConfigu}) admits a unimodular triangulation $T$.\footnote{If $T$ is not unimodular, the description of the combinatorial structure of $\sol_{M_A(\delta),z}$ and that of $\Homo_n(V_z^{an};\mathcal{L})$ are more complicated (\cite{GM}).}
Let us introduce notation. For any subset $\tau\subset\{1,\dots,N\}$, we set $z_{\tau}:=(z_i)_{i\in\tau}$. We write $\bs$ for the complement $\{ 1,\dots,N\}\setminus \s$. For any vector $v=(v_1,v_2,\dots)$ and a univariate scalar-valued function $F$, we define $F(v)$ to be the product of values $F(v):=F(v_1)F(v_2)\cdots$. With this notation, for any $(n+k)$-simplex $\s$, we set
\begin{equation}\label{seriesphi}
\varphi_{\s}(z;\delta):=
z_\sigma^{-A_\sigma^{-1}\delta}
\sum_{{\bf m}\in\Z_{\geq 0}^{\bs}}\frac{z_\sigma^{-A_\sigma^{-1}A_{\bar{\sigma}}{\bf m}}z_{\bar{\sigma}}^{\bf m}}{\Gamma({\bf 1}-A_\sigma^{-1}(\delta+A_{\bar{\sigma}}{\bf m})){\bf m!}}
\end{equation}
where ${\bf 1}$ is a vector of length $n+k$ whose entries are all unity. By a direct computation, we can show that $\varphi_{\s}(z;\delta)$ is a solution of $M_A(\delta)$. We say that $\delta$ is very generic if any entry of the vector $A_\sigma^{-1}(\delta+A_{\bar{\sigma}}{\bf m})$ is non-integral. It is easy to see that if $\delta$ is very generic, $\delta$ must be non-resonant.  \begin{prop}
If $(-\log|z_1|,\dots,-\log|z_N|)$ is in a sufficiently far translation of the cone $C_T$ inside itself, the series (\ref{seriesphi}) is convergent. Moreover, if $\delta$ is very generic, $\{\varphi_{\s}(z;\delta)\}_{\s\in T}$ is a basis of $\sol_{M_A(\delta),z}$.
\end{prop}

Through the isomorphism $\Homo_n(V_z^{an};\mathcal{L})\simeq\sol_{M_A(\delta),z}$, we have the basis $\{[\Gamma_\s]\}_{\s\in T}$ of $\Homo_n(V_z;\mathcal{L})$ corresponding to $\{\varphi_{\s}(z;\delta)\}_{\s\in T}$. Similarly, we have a dual basis $\{[\Gamma_\s^\vee]\}_{\s\in T}$ of $\Homo_n(V_z^{an};\mathcal{L}^\vee)$. An important point here is that these bases are orthogonal bases with respect to the homology intersection form. Namely, we have the orthogonality relation $\langle [\Gamma_{\s_1}],[\Gamma_{\s_2}^\vee]\rangle_h=0$ if $\s_1\neq\s_2$. The remaining homology intersection number $\langle [\Gamma_{\s}],[\Gamma_{\s}^\vee]\rangle_h$ is also explicitly given by
\begin{equation}\label{eqn:SigmaIntersection}
\langle [\Gamma_{\s}],[\Gamma_{\s}^\vee]\rangle_h=C(\gamma;\sigma)\displaystyle\prod_{l:|\s^{(l)}|>1}\left\{(1-e^{2\pi\ii\gamma_l})\prod_{i\in\s^{(l)}}\left(1-e^{-2\pi\ii{}^t{\bf e}_iA_\s^{-1}\delta}\right)\right\},
\end{equation}
where $C(\gamma;\s)$ is a constant depending only on $\gamma_1,\dots,\gamma_k$ and $\s$ (\cite[Theorem 7.5]{MH}).

For any complex vectors ${\bf v}=(v_1,v_2,\dots)$ and ${\bf w}=(w_1,w_2,\dots)$ of equal length, we set $({\bf v})_{\bf w}:=\frac{\Gamma({\bf v+w})}{\Gamma({\bf v})}$. Finally, we set $h^{\bf b}:=h_1^{b_1}\cdots h_k^{b_k}$ for any element ${\bf b}=(b_1,\dots,b_k)\in\Z^k$. In view of Riemann-Hodge bilinear relation (\ref{eqn:expansion}), we obtain an expansion theorem of the cohomology intersection number.

\begin{theorem}[Theorem 8.1 of \cite{MH}]\label{thm:QuadraticRelation}
Suppose that four vectors ${\bf a},{\bf a}^\prime\in\Z^{n},{\bf b},{\bf b}^\prime\in\Z^{k}$ and a unimodular regular triangulation T\footnote{As for the case when $T$ is not unimodular, see \cite[Theorem 2.6]{GM}.} are given. If the parameter $\delta$ is generic so that $\gamma_l\notin\Z$ for any $l=1,\dots,k$ and $\delta$, $
\begin{pmatrix}
\gamma-{\bf b}\\
c+{\bf a}
\end{pmatrix}
$  and 
$
\begin{pmatrix}
\gamma+{\bf b}^\prime\\
c-{\bf a}^\prime
\end{pmatrix}
$ are very generic, then, one has an identity
\begin{align}
&(-1)^{|{\bf b}|+|{\bf b}^\prime|}\gamma_1\cdots\gamma_k(\gamma-{\bf b})_{\bf b}(-\gamma-{\bf b}^\prime)_{{\bf b}^\prime}\times\nonumber\\
&\sum_{\s\in T}\frac{\pi^{n+k}}{\sin\pi A_\s^{-1}\delta}\varphi_{\s}\left(z;
\begin{pmatrix}
\gamma-{\bf b}\\
c+{\bf a}
\end{pmatrix}
\right)\varphi_{\s}\left(z;
\begin{pmatrix}
-\gamma-{\bf b}^\prime\\
-c+{\bf a}^\prime
\end{pmatrix}
\right)\nonumber\\
=&\frac{\langle x^{\bf a}h^{\bf b}\frac{dx}{x},x^{{\bf a}^\prime}h^{{\bf b}^\prime}\frac{dx}{x}\rangle_{ch}}{(2\pi\ii)^n}\label{eqn:RCIN}
\end{align}
for any $z$ such that $(-\log|z_1|,\dots,-\log|z_N|)$ is in a sufficiently far translation of the cone $C_T$ inside itself.
\end{theorem}

\noindent
Since cohomology classes $[x^{\bf a}h^{\bf b}\frac{dx}{x}]$ generate the algebraic de Rham cohomology group $\Homo_{dR}^n(V_z;\nabla_x)$, Theorem \ref{thm:QuadraticRelation} gives a closed formula of any cohomology intersection number. However, Theorem \ref{main} implies that the cohomology intersection number is a priori a rational function while the formula (\ref{eqn:RCIN}) is, in general, an infinite series.

Let us illustrate how Theorem \ref{thm:QuadraticRelation} is used to evaluate cohomology intersection number. We take free bases $\{[\omega_i(z)]\}_{i=1}^r\subset\HdR^n$ and $\{[\omega_i^\vee(z)]\}_{i=1}^r\subset\HdR^{n\vee}$. Suppose that $I$ is a non-zero rational solution of the secondary equation (\ref{SecondaryEq2}). In view of Theorem \ref{main}, there is a complex constant $C$ such that $I_{ch}=\left(\langle [\omega_i(z)],[\omega_j^\vee(z)]\rangle_{ch}\right)_{i,j=1}^r=C\cdot I$. We choose some $i$ and $j$ and focus on $\langle [\omega_i(z)],[\omega_j^\vee(z)]\rangle_{ch}$. If $I_{ij}(z)$ denotes the $(i,j)$-entry of $I$, we have 
\begin{equation}\label{eqn:Triv}
\langle [\omega_i(z)],[\omega_j^\vee(z)]\rangle_{ch}=C\cdot I_{ij}(z).
\end{equation}
We may assume that $I_{ij}(z)$ is a non-zero function. Since the cohomology classes $[\omega_i(z)]$ and $[\omega_j^\vee(z)]$ are expanded into $\C(z)$-linear combination of the cohomology classes of the form $[x^{\bf a}h^{\bf b}\frac{dx}{x}]$, we can apply the formula (\ref{eqn:RCIN}) to obtain a Laurent expansion of $\langle [\omega_i(z)],[\omega_j^\vee(z)]\rangle_{ch}$. Then, we substitute a particular value $z=z_0$ in (\ref{eqn:Triv}), which determines the constant $C$. The point $z_0$ is usually taken as the ``center'' of the Laurent expansion. When the value $I_{ij}(z_0)$ diverges or vanishes, we divide or multiply (\ref{eqn:Triv}) by a suitable polynomial factor in $z$ before substitution.

We conclude this subsection by citing a theorem on an arithmetic property of the cohomology intersection number. We define a field $\Q(\delta)$ as a  field extension $\Q(\delta):=\Q(\gamma_1,\dots,\gamma_k,c_1,\dots,c_n)$ of $\Q$.
\begin{theorem}[Theorem 2.9 of \cite{GM} and Theorem 3.5 of \cite{MT}]
Suppose that $\delta$ is non-resonant and $\gamma_l\notin\Z$. Then, for any $P_1,P_2\in\mathbb{Q}(\delta)\langle z,\partial_z\rangle$, the cohomology intersection number $\frac{\langle P_1\cdot\frac{dx}{x},P_2\cdot\frac{dx}{x}\rangle_{ch}}{(2\pi\ii)^n}$ belongs to the field $\mathbb{Q}(\delta)(z)$. 
\end{theorem}

\noindent
The theorem above guarantees that we do not need any field extension of $\mathbb{Q}(\delta)$ when we compute the cohomology intersection number.

\section{GKZ system and algorithms}

In this section, we set $\beta:=-\delta$. With this notation, we put $H_A(\beta):=M_A(\delta)$. This is because we use some results from \cite{dojo-en} and \cite{SSTip} where the hypergeometric ideal is denoted by $H_A(\beta)$ while it is denoted by $M_A(\delta)$ in our main references \cite{MH}, \cite{MT}.

\subsection{An algorithm of computing connection matrices}\label{subsec:5.1}

Let $\omega_{q}$ be the differential form
\begin{equation} \label{eq:basis_format}
  \prod_{l=1}^k h_l^{-q'_l} x^{q''} \frac{dx}{x}, \quad
  q=(q',q'') \in \Z^{k} \times \Z^n .
\end{equation}
In view of Theorem \ref{thm:thm2}, there exists a basis of the twisted cohomology group of which elements are of the form $\omega_{q}$
when $\delta$ is non-resonant and $\gamma_l\notin\Z$. Such a basis is even algorithmically computable.
Let $\{ [\omega_{q}] \,|\, q \in Q \}$ be a basis of the twisted cohomology group.
We set $\pd{i} = \frac{\partial}{\partial z_i}$.
We will give an algorithm to find the connection matrix 
$ \nabla^{GM}_{\pd{i}} \omega = \omega \Omega_i$
with respect to this basis $\omega = ( [\omega_{q_1}],\dots, [\omega_{q_r}])$ where $Q=\{ q_1,\dots,q_r\}$. In the theory of differential equations, it is more common to consider Pfaffian matrix $P_i:={}^t\Omega_i$ instead of connection matrix. Note that algorithms to translate a given holonomic ideal
to a Pfaffian system are well known
(see, e.g., \cite[Chap 6]{dojo-en}).
In the following, we explain how we compute the matrix $P_i$ by means of computer algebra, which was proposed in \cite{MT2}.

The main point of our method lies in the use of the following contiguity relation
\begin{equation} \label{eq:annihilating}
\frac{1}{{\bf a}'_i \cdot (\beta-q)} \pd{i} \cdot [\omega_{q}] =[ \omega_{q'} ], \quad
 q' = q + {\bf a}(i)
\end{equation}
where ${\bf a}'_i$ is the column vector that the first $k$ elements are 
equal to those of ${\bf a}(i)$ and the last $n$ elements are $0$.
For example, ${\bf a}'_1 = {}^t(1, 0, \ldots, 0)$,
${\bf a}'_2 = {}^t(1, 0, \ldots, 0)$, $\ldots$,
${\bf a}'_{N_1+1} = {}^t(0,1,0, \ldots,0)$, $\ldots$,
${\bf a}'_N = {}^t(0,\ldots,0,\overset{k-th}{1},0,\ldots,0)$.

In \cite[Algorithm 3.2]{SSTip}, an algorithm to obtain 
an operator $C_i$ satisfying
\begin{equation}
  C_i \pd{i} - b_i(\beta) = 0 \quad {\rm mod}\ H_A(\beta)
\end{equation}
is given.
The polynomial $b_i$ is a $b$-function in the direction $i$ \cite[Th 3.2]{SSTip}.
Note that the algorithm outputs the operator $C_i$ in
$\C \langle z_1, \ldots, z_N, \pd{1}, \ldots, \pd{N} \rangle$,
which does not depend on the parameter $\beta$.
We have the following inverse contiguity relation
\begin{equation} \label{eq:creating}
 \frac{{\bf a}_i' \cdot (\beta-q'')}{b_i(\beta-q'')}C_i \cdot [\omega_{q}] = [ \omega_{q''}], \quad
 q'' = q - {\bf a}(i) .
\end{equation}

\begin{example}\rm  \label{example:gauss} (Gauss hypergeometric function ${}_2F_1$.)
Put
\begin{equation}
A = \left(
\begin{array}{cc|cc}
1 & 1 & 0 & 0 \\ \hline
0 & 0 & 1 & 1 \\ \hline
0 & 1 & 0 & 1 \\
\end{array}
\right).
\end{equation}
The integral (\ref{eqn:EulerInt}) in question takes the form
\begin{equation}
\int_\Gamma h_1^{-\gamma_1}h_2^{-\gamma_2}x^c\omega
\end{equation}
where $h_1 = z_1 + z_2 x$ and $h_2 = z_3 + z_4 x$. We can show that
$\{ [\omega_{(1,0,0)}], [\omega_{(0,1,0)}] \}$
is a basis of the de Rham cohomolgy group $\Homo^1_{dR}(V_z;\nabla_x)$.
This $A$ is normal and the $b$-function $b_4(s) \in \Q[s_1,s_2,s_3]$ 
for the direction $z_4$ is
$ b_4(s) = s_2 s_3$.
Then, $C_4 =  z_2 z_3 \pd{1} + (\theta_2 + \theta_3 + \theta_4)z_4$
where $\theta_i = z_i \pd{i}$
by reducing $(\theta_3+\theta_4)(\theta_2+\theta_4)$
by the toric ideal $I_A=\langle \underline{\pd{2}\pd{3}}-\pd{1}\pd{4} \rangle$ 
(see Algorithm 3.2 of \cite{SSTip}). 
\end{example}

Our algorithm to find a Pfaffian system with respect to a given
basis of the twisted cohomology group is as follows. 
\begin{algorithm} \label{alg:pf} \rm \quad \\
Input: $\{ [\omega_q] \,|\, q \in Q\}$, a basis of the twisted cohomology group. 
A direction (index) $i$. \\
Output: $P_i$, the coefficient matrix of the Pfaffian system. 
\begin{enumerate}
\item Compute a Gr\"obner basis $G$ of $H_A(\beta)$ in the ring of differential operators
with rational function coefficients.
Let $S$ be a column vector of the standard monomials with respect to $G$.
\item Put
\begin{equation} \label{eq:composite}
 F(Q) = {}^t(F(q) \,|\, q \in Q), \quad
  F(q) =  \prod_{r_i < 0} C_i^{-r_i} \prod_{r_i > 0} \pd{i}^{r_i} \frac{1}{BB'}, \quad
 q = \sum_{i=1}^{N} r_i {\bf a}(i)
\end{equation}
It is a vector with entries in the ring of differential operators
and the order of the product is $ i = N, N-1, \ldots, 3, 2, 1$.
In other words, we apply operators from $\pd{1}$.
The polynomial $B$ is derived from the coefficient of the contiguity relation
(\ref{eq:creating})
and is equal to
\begin{eqnarray} \label{eq:B_i}
B   &=& \prod_{j=1, r_j<0}^{N}  \frac{b_j(\beta_j'+{\bf a}(j))}
                                 {{\bf a}'_j \cdot (\beta_j'+{\bf a}(j))}
                            \frac{ b_j(\beta_j'+2{\bf a}(j))}
                                 {{\bf a}'_j \cdot (\beta_j'+2{\bf a}(j))} \cdots 
                            \frac{b_j(\beta_j'+(-r_j){\bf a}(j))}
                                 {{\bf a}'_j \cdot (\beta_j'+(-r_j){\bf a}(j))}, \\
\beta_j' &=& \beta-\sum_{r_l>0} r_l {\bf a}(l) + \sum_{l=1, r_l<0}^{j-1} (-r_l) {\bf a}(l).
\end{eqnarray}
The polynomial $B'$ comes from the denominator of the contiguity relation
(\ref{eq:annihilating}) and is equal to
\begin{eqnarray} \label{eq:B'_i}
B' &=& \prod_{j=1, r_j>0}^{N}  \left( {\bf a}'_j \cdot (\beta_j') \right)
   \left( {\bf a}'_j \cdot(\beta_j'-{\bf a}(j))\right) \cdots 
   \left( {\bf a}'_j \cdot (\beta_j'-(r_j-1) {\bf a}(j)) \right), \\
\beta_j' &=& \beta-\sum_{r_l>0, l < j} r_l {\bf a}(l) .
\end{eqnarray}
\item Compute the normal form of the vectors
$\pd{i} F(Q)$ and $F(Q)$.
Write the normal forms of them as
$P'S$ and $P''S$ respectively
where $P'$ and $P''$ are matrices with rational function entries.
\item Output $P_i =  P'(P'')^{-1}$.
\end{enumerate}
\end{algorithm}
The matrix $P''$ is invertible if and only if 
the given set of cohomology classes $\{ [\omega_q] \}$ 
is a basis of the twisted cohomology group.

\begin{example}\label{example:gauss2}
\rm
This is a continuation of Example \ref{example:gauss}.
We have
${}^t(1,0,0) = {\bf a}(1)$ and ${}^t(0,1,0) = {\bf a}(3)$.
Then, the basis of the twisted cohomology group $F(Q)$
is expressed as
$ F(Q) = {}^t(\pd{1}/\beta_1,\pd{3}/\beta_2)$
and
$\pd{4} F(Q) = {}^t(\pd{4}\pd{1}/\beta_1,\pd{4}\pd{3}/\beta_2)$.
We can obtain a Gr\"obner basis whose set of the standard monomials is
$\{ \pd{4}, 1 \}$
by the graded reverse lexicographic order such that $\pd{i} > \pd{i+1}$.
We multiply $\beta_1 \beta_2$ to $F(Q)$ and $\pd{4} F(Q)$ in order to avoid
rational polynomial arithmetic.
Then, the normal form, for example, of
$\beta_2 \pd{1}$ is 
\newline
$\frac{1}{z_1 z_4-z_2 z_3} \left((\beta_1(\beta_1+\beta_2) z_4)\pd{4}-\beta_2^2 \beta_3\right)$.
By computing the other normal forms, we obtain
the matrix
\begin{equation}
 P_4=
\left(
\begin{array}{cc}
\frac{  {\beta}_{2}  {z}_{1}}{   {z}_{1}  {z}_{4}-  {z}_{2}  {z}_{3}}& \frac{ -  {\beta}_{2}  {z}_{3}}{   {z}_{1}  {z}_{4}-  {z}_{2}  {z}_{3}} \\
\frac{ -   {\beta}_{1}  {z}_{1}  {z}_{2}}{   {z}_{1}   {z}_{4}^{ 2} -   {z}_{2}  {z}_{3}  {z}_{4}}& \frac{    {\beta}_{3}  {z}_{1}  {z}_{4}+   (  {\beta}_{1}- {\beta}_{3})  {z}_{2}  {z}_{3}}{   {z}_{1}   {z}_{4}^{ 2} -   {z}_{2}  {z}_{3}  {z}_{4}}
\end{array}
\right).
\end{equation}
\end{example}

\begin{example} \rm \label{example:3F2}
(${}_3F_2$, see, e.g., \cite[p.224]{SST}, \cite{Ohara-Sugiki-Takayama}.)
Let 
$A = 
\left(
\begin{array}{cc|cc|cc}
 1&  1& 0& 0& 0& 0 \\ \hline
0& 0&  1&  1& 0& 0 \\ \hline
0& 0& 0& 0&  1&  1 \\ \hline
 1& 0& 0&  1& 0& 0 \\
0& 0&  1& 0& 0&  1 \\
\end{array}
\right)
$.
The integral (\ref{eqn:EulerInt}) in question takes the form 
\begin{equation}
 \int_\Gamma (z_1 x_1 + z_2)^{-\gamma_1} (z_3 x_2 + z_4 x_1)^{-\gamma_2}
             (z_5 + z_6 x_2)^{-\gamma_3} x_1^{c_1} x_2^{c_2} \omega.
\end{equation}
We set
\begin{equation}
\omega_1 = \frac{dx_1 dx_2}{(z_1 x_1 + z_2)x_1 x_2},
  \omega_2 = \frac{dx_1 dx_2}{(z_5  + z_6 x_2)x_1 x_2},
  \omega_3 = \frac{dx_1 dx_2}{(z_3 x_2 + z_4 x_1)x_1 x_2}.
\end{equation}
It can be verified that $\{ [\omega_1],[\omega_2],[\omega_3]\}$ is a basis of the de Rham cohomology group $\Homo^2(V_z;\nabla_x)$. When $z_2=-1, z_3=z_4=z_5=z_6=1$, the coefficient matrix for $z_1$
for the basis 
$\{ [\omega_1],[\omega_2],[\omega_3]\}$ 
is 
\begin{equation}
 P_1=
\left(
\begin{array}{ccc}
\frac{   \beta_{4}  z_{1}+  \beta_{2}+  \beta_{3}  - \beta_{4}- \beta_{5}}{   z_{1}(z_1 - 1)}& \frac{   \beta_{3}  (\beta_{1}+   \beta_{2}-  \beta_{4} )}{   \beta_{1}   z_{1}( z_{1}-1)}& \frac{   \beta_{2}(\beta_2   - \beta_{4}  - \beta_{5}- 1)  }{   \beta_{1}   z_{1}( z_{1}-1)} \\
\frac{  (  \beta_{2}+  \beta_{3}- \beta_{5})  \beta_{1}}{   \beta_{3}  (z_{1}-1)}& \frac{   \beta_{1}  z_{1}+  \beta_{2}- \beta_{4}}{   z_{1}( z_{1}-1)}& \frac{   \beta_{2}(\beta_2  - \beta_{4}  - \beta_{5}- 1) }{   \beta_{3}   z_{1}( z_{1}-1)} \\
\frac{  (  - \beta_{2}  - \beta_{3}+ \beta_{5})  \beta_{1}}{   \beta_{2} ( z_{1}-1)}& \frac{ \beta_3(\beta_4 -   \beta_{1}  -   \beta_{2})}{   \beta_{2}  (z_{1}-1)}& \frac{  - \beta_{2}+  \beta_{4}+  \beta_{5}+ 1}{  z_{1}- 1}
\end{array}
\right)
\end{equation}
The result can be obtained in a few seconds.
\end{example}

\subsection{An algorithm of finding the cohomology intersection matrix} \label{sec:5}

\begin{theorem} {\rm \cite{MT}}
Given 
a matrix $A=(a_{ij})$ as in (\ref{CayleyConfigu}).\footnote{In \cite{MT}, the matrix $A$ is assumed to have a unimodular regular triangulation. This technical assumption is not necessary in view of \cite[Theorem 2.6]{GM}}
When parameters are non-resonant and $\gamma_l\notin\Z$, 
the intersection matrix of
the twisted cohomology group of the GKZ system associated to the matrix $A$
can be algorithmically determined.
\end{theorem}

We write $\Omega_i$ for the coefficient matrix of $\Omega$
with respect to the $1$-form $dz_i$.
The algorithm we propose is summarized as follows.

\begin{algorithm} {\rm (A modified version of the algorithm in \cite{MT}.)}
 
Input: Free bases $\{[\phi_j]\}_j\subset\HdR^n\restriction_U$,
$\{[\psi_j]\}_j\subset\HdR^{n\vee}\restriction_U$ which are expressed as (\ref{eq:basis_format}).

Output: The secondary equation (\ref{SecondaryEq2}) and the
cohomology intersection matrix $I_{ch}=(\langle[\phi_i],[\psi_j]\rangle_{ch})_{i,j}$.

\begin{enumerate}
\item Obtain a Pfaffian system with respect to the given bases $\{[\phi_j]\}_j$ and $\{[\psi_j]\}_j$, i.e.,
obtain matrices $\Omega_i=(\omega_{ijk})$ and $\Omega_i^\vee=(\omega^\vee_{ijk})$ so that the equalities
\begin{equation}
 \pd{i}[\phi_j]=\sum_k\omega_{ikj}[\phi_k],\ \ \ \pd{i}[\psi_j]=\sum_k\omega^\vee_{ikj}[\psi_k]
\end{equation}
hold by Algorithm \ref{alg:pf}.
\item Find a non-zero rational function solution $I$ of the secondary equation
\begin{equation}
 \pd{i} I - {}^t\Omega_i I - I \Omega_i^\vee =0,
\quad i=1, \ldots, N. 
\end{equation}
To be more precise, see, e.g., \cite{Barkatou}, \cite{IC-proj}, \cite{Oaku-Takayama-Tsai} and references therein.
\item Determine the scalar multiple of $I$ by Theorem \ref{thm:QuadraticRelation} or by \cite[Theorem 2.6]{GM}.
\end{enumerate}

\end{algorithm}

\begin{example}
\rm This is a continuation of Example \ref{example:gauss} and Example \ref{example:gauss2}. In this case, we set $\omega_1=\omega_1^\vee=\omega_{(1,0,0)}$ and $\omega_2=\omega_2^\vee=\omega_{(0,1,0)}$. By solving the secondary equation (for example, using \cite{IC-proj}), we can verify that $I_{ch}=(\langle [\omega_i],[\omega^\vee_j]\rangle_{ch})_{i,j=1}^2$ is a constant matrix when $z_1=z_2=z_3=1$. Therefore, we can obtain the exact values of these entries by taking a unimodular regular triangulation $T=\{ 123,234\}$ and substituting $z_4=0$ in Theorem \ref{thm:QuadraticRelation}. Thus, we get a correct normalization of $I_{ch}$ and the matrix $\frac{I_{ch}|_{z_1=z_2=z_3=1}}{2\pi\ii}$ is given by
\begin{equation}
\begin{pmatrix}
\frac{1}{\beta_1}-\frac{1}{\beta_3}&-\frac{1}{\beta_3}\\
-\frac{1}{\beta_3}&\frac{1}{\beta_2}-\frac{1}{\beta_3}
\end{pmatrix}.
\end{equation}
\end{example}

\begin{example}
\rm This is a continuation of Example \ref{example:3F2}. We want to evaluate the cohomology intersection matrix $I_{ch}=(\langle [\omega_i],[\omega_j]\rangle_{ch})_{i,j=1}^3$. By solving the secondary equation, we can verify that $(1,1),$ $(1,2),$ $(2,1),$ $(2,2)$ entries of $I_{ch}|_{-z_2=z_3=z_4=z_5=z_6=1}$ are all independent of $z_1$. Therefore, we can obtain the exact values of these entries by taking a unimodular regular triangulation $T=\{ 23456,12456,12346\}$ and substituting $z_1=0$ in Theorem \ref{thm:QuadraticRelation}. Thus, the matrix $\frac{I_{ch}|_{-z_2=z_3=z_4=z_5=z_6=1}}{(2\pi\ii)^2}$ is given by

\tiny
\begin{equation}
\left[ \begin {array}{ccc} r_{11}
&{\frac 
{\beta_4+\beta_5}{ \left( \beta_2-\beta_4-\beta_5 \right) \beta_5\,\beta_4}}&{\frac {\beta_1\,\beta_4\,z_1+\beta_2\,\beta_4
\,z_1-{\beta_4}^{2}z_1-\beta_4\,\beta_5\,z_1-\beta_5\beta_3}{ \left( \beta_2-\beta_4-\beta_5+1 \right)  \left( \beta_2-\beta_4-\beta_5 \right) \beta_5\,\beta_4}}
\\
\noalign{\medskip}{\frac {\beta_4+\beta_5}{ \left( \beta_2-\beta_4-\beta_5 \right) \beta_5\,\beta_4}}&r_{22}&-{\frac {\beta_1\,\beta_4\,z_1-\beta_5\,\beta_2-\beta_5\,
\beta_3+\beta_5\,\beta_4+{\beta_5}^{2}}{ \left( \beta_2-\beta_4-
\beta_5+1 \right)  \left( \beta_2-\beta_4-\beta_5 \right) \beta_5
\,\beta_4}}\\
\noalign{\medskip}{\frac {\beta_1\,\beta_4\,z_1+
\beta_2\,\beta_4\,z_1-{\beta_4}^{2}z_1-\beta_4\,\beta_5
\,z_1-\beta_5\,\beta_3}{ \left( \beta_2-\beta_4-\beta_5-1
 \right)  \left( \beta_2-\beta_4-\beta_5 \right) \beta_5\,\beta_4
}}&-{\frac {\beta_1\,\beta_4\,z_1-\beta_5\,\beta_2-\beta_5\,
\beta_3+\beta_5\,\beta_4+{\beta_5}^{2}}{ \left( \beta_2-\beta_4-
\beta_5-1 \right)  \left( \beta_2-\beta_4-\beta_5 \right) \beta_5
\,\beta_4}}&r_{33}\end {array} \right]
\end{equation}
\normalsize
where
\begin{equation}
r_{11} =  -\frac{      (   {\beta}_{4}  {\beta}_{2}+  (  {\beta}_{4}+ {\beta}_{5})  {\beta}_{3})  {\beta}_{1}+  {\beta}_{4}   {\beta}_{2}^{ 2} +  (    {\beta}_{4}  {\beta}_{3}-  {\beta}_{4}^{ 2} -  {\beta}_{5}  {\beta}_{4})  {\beta}_{2}+  (  -  {\beta}_{4}^{ 2} -  {\beta}_{5}  {\beta}_{4})  {\beta}_{3}} {     {\beta}_{5}  {\beta}_{4}   {\beta}_{1}(   {\beta}_{2}- {\beta}_{4}- {\beta}_{5})  (   {\beta}_{2}+ {\beta}_{3}- {\beta}_{5})  }
\end{equation}
\begin{equation}
r_{22}=-{\frac {\beta_1\,\beta_2
\,\beta_5+\beta_1\,\beta_3\,\beta_4+\beta_1\,\beta_3\,\beta_5-\beta_1\,\beta_4\,\beta_5-\beta_1\,{\beta_5}^{2}+\beta_5\,{\beta_2}^{2}+\beta_2\,\beta_3\,\beta_5-\beta_2\,\beta_4\,\beta_5-\beta_2\,{\beta_5}^{2}}{\beta_3\, \left( \beta_1+\beta_2-\beta_4
 \right)  \left( \beta_2-\beta_4-\beta_5 \right) \beta_5\,\beta_4
}}
\end{equation}
\begin{equation}
r_{33}=-\frac {\alpha_0z_1^2-2\,\beta_1\,\beta_3\,\beta_4\,\beta_5\,
z_1+\alpha_2}{ \left( \beta_2-\beta_4-\beta_5-1 \right)  \left( \beta_2-\beta_4-\beta_5+1 \right) \beta_2\, \left( \beta_2-\beta_4-\beta_5 \right) \beta_5\,\beta_4
}
\end{equation}
\begin{equation}
\alpha_0={\beta_1}^{2}\beta_2\,\beta_4
-{\beta_1}^{2}\beta_4\,\beta_5+\beta_1\,{\beta_2}
^{2}\beta_4-\beta_1\,\beta_2\,{\beta_4}^{2}-2\,\beta_1\,\beta_2\,\beta_4\,\beta_5+\beta_1\,{\beta_4}^{2}\beta_5+\beta_1\,\beta_4\,{\beta_5}^{2}
\end{equation}
\begin{equation}
\alpha_2={\beta_2}^{2}\beta_3\,\beta_5+\beta_2\,{\beta_3}^{2}\beta_5-2\,\beta_5\,\beta_4\,\beta_3\,\beta_2-\beta_2\,\beta_3\,
{\beta_5}^{2}-{\beta_3}^{2}\beta_4\,\beta_5+\beta_3\,{\beta_4}^{
2}\beta_5+\beta_3\,\beta_4\,{\beta_5}^{2}
\end{equation}
\end{example}

\section{$L^2$-cohomology intersection pairing and an integral of a product of powers of absolute values of polynomials} \label{sec:3}

\subsection{$L^2$-cohomology intersection pairing}\label{sec:61}

We use the same notation as \S\ref{sec:Basic}. We want to understand an integral

\begin{equation}\label{eqn:AbsInt}
I(\alpha)=\int_{\C^n}|\Phi|^{2}\omega\wedge\bar{\eta}
\end{equation}
as a meromorphic function of $\alpha$ for some $\omega,\eta\in\Omega^n(V^{alg})$. In order to analyze (\ref{eqn:AbsInt}), we employ the language of $L^2$-cohomology groups. We assume that $\alpha_1,\dots,\alpha_k\in\R$. We first remark that $\mathcal{L}$ is trivially a variation of Hodge structure of weight $0$ (\cite{KKPoincare}). Moreover, there is a polarization given by $\mathcal{L}\otimes \overline{\mathcal{L}}\ni a\Phi\otimes \overline{b\Phi}\mapsto a\bar{b}\in\C$.  Here, the symbol $\overline{\mathcal{L}}$ denotes the complex conjugate of the local system $\mathcal{L}$. We consider a smooth  projective compactification $X$of $V$ so that the complement $D:=X\setminus V$ is a normal crossing divisor. Let us fix a K\"ahler metric $g$ on $X$ which is asymptotically equivalent to Poincar\'e metric near the boundary $D$. Namely, our K\"ahler metric $g$ dominates and is dominated by a positive multiple of 
\begin{equation}\label{eqn:Kahler}
\sum_{j\leq l}\frac{dx_jd\bar{x}_j}{(|x_j|\log|x_j|)^2}+\sum_{j>l}dx_jd\bar{x}_j
\end{equation}
on the coordinate system $x$ near the boundary such that $D=\{ x_1\cdots x_l=0\}$. Note that the volume form induced from (\ref{eqn:Kahler}) is a constant multiple of $\wedge_{j\leq l}|x_j|^{-1}(\log|x_j|)^{-1}dx_j\wedge d\bar{x}_j\wedge\wedge_{j> l}dx_j\wedge d\bar{x}_j$. For a polarized variation of Hodge structures $H$ on $V$, the symbol $\mathcal{L}^p_{(2)}(H)$ denotes the sheaf of $L^2$-differential $p$-forms with values in $H$ (\cite[Definition 5.3.1]{KKPoincare}). We cite the result of \cite[Theorem 5.4.1]{KKPoincare}.
\begin{theorem}[\cite{KKPoincare}]
The complex $(\mathcal{L}_{(2)}^\bullet(H),d)$ is quasi-isomorphic to the minimal extension ${}^\pi H$ of $H$ on $X$.
\end{theorem}
We consider a variation of Hodge structures $H=\mathcal{L}^\vee$ and set
\begin{equation}
\Homo_{(2)}^n(V^{an},\mathcal{L}^\vee):=\mathbb{H}^n(X^{an};(\mathcal{L}_{(2)}^\bullet(H),d)).
\end{equation}
We describe the $L^2$-intersection pairing
\begin{equation}\label{eqn:L2Pairing}
\langle\bullet,\bullet\rangle:\Homo^n_{(2)}(V^{an},\mathcal{L}^\vee)\times \Homo^n_{(2)}(V^{an},\overline{\mathcal{L}^\vee})\rightarrow\C,
\end{equation}
which was given in \cite[Theorem 6.4.2]{KKPoincare}. If we use the resolution $(\mathcal{L}_{(2)}^\bullet(H),d)$ of ${}^\pi H$, (\ref{eqn:L2Pairing}) is induced from the local duality pairing $\mathcal{L}_{(2)}^n(H)\otimes \mathcal{L}_{(2)}^n(\bar{H})\ni (\xi\otimes\Phi^{-1})\otimes(\eta\otimes \bar{\Phi}^{-1})\mapsto \xi\wedge\eta\in \Db^{2n}_{X^{an}}$ where $\Db^{2n}_{X^{an}}$ is the sheaf of $2n$-currents on $X^{an}$. For our purpose, it is more convenient to use another resolution. We write $\Db^{{\rm mod}D,p}_{X^{an}}$ for the sheaf of $p$-currents with moderate growth along $D$. Let us consider the quasi-isomorphism $(\Db_{X^{an}}^{{\rm mod}D,\bullet},\nabla_x)\rightarrow(\Db_{X^{an}}^{{\rm mod}D,\bullet}\otimes H,d)$ given by the correspondence $\Db_{X^{an}}^{{\rm mod}D,p}\ni\varphi\mapsto \varphi\Phi\otimes\Phi^{-1}\in\Db_{X^{an}}^{{\rm mod}D,p}\otimes H$. This morphism does not depend on a particular choice of a branch of $\Phi$ and therefore it is well-defined. Since $\mathcal{L}_{(2)}^p(H)$ is a subsheaf of $\Db_{X^{an}}^{{\rm mod}D,p}\otimes H$, we can define a subsheaf $\mathcal{L}_{(2)}^p$ of $\Db_{X^{an}}^{{\rm mod}D,p}$ so that there is a quasi-isomorphism $(\mathcal{L}_{(2)}^{\bullet},\nabla_x)\rightarrow(\mathcal{L}_{(2)}(H),d)$. That a measurable $p$-form $\varphi$ be a section of $\mathcal{L}_{(2)}^p$ is characterized by the condition that both $|\Phi|^2\varphi\wedge*\overline{\varphi}$ and $|\Phi|^2\nabla_x\varphi\wedge*\overline{\nabla_x\varphi}$ are integrable. Here, $*$ is the Hodge star operator. In sum, we obtain an identity
\begin{equation}
\Homo^n_{(2)}(V^{an},\mathcal{L}^\vee)=\mathbb{H}^n(X^{an}; (\mathcal{L}_{(2)}^{\bullet},\nabla_x)).
\end{equation}
We can describe the pairing (\ref{eqn:L2Pairing}) by the formula $\Homo^n_{(2)}(X^{an},\mathcal{L}^\vee)\times \Homo^n_{(2)}(X^{an},\overline{\mathcal{L}^\vee})\ni [\omega]\otimes[\bar{\eta}]\mapsto\int_{X^{an}}|\Phi|^2\omega\wedge \bar{\eta}\in\C$.

Now we focus on the case when the regularization condition is satisfied. Namely, we assume that the canonical morphisms $j_!H\rightarrow{}^\pi H$ and ${}^\pi H\rightarrow\R j_*H$ are isomorphisms. The regularization condition is again a generic condition on parameters $\alpha$. Under this condition, we have the following commutative diagram

\begin{equation}\label{eqn:CD}
\xymatrix{
\Homo^n_{dR,c}(V^{an};\nabla_x) \ar[d]_{{\rm can}_{(2)}} \ar[r]^{{\rm can}}&\Homo^n_{dR}(V^{an};\nabla_x) \\
\Homo_{(2)}^n(V^{an};\mathcal{L}^\vee) \ar[ur]_{\iota} & 
}.
\end{equation}
Here, the morphisms ${\rm can}$ and ${\rm can}_{(2)}$ are induced from the canonical morphisms $j_!H\rightarrow\R j_* H$ and $j_!H\rightarrow{}^\pi H$ respectively, and the morphism $\iota:\Homo_{(2)}^n(V^{an};\mathcal{L}^\vee)\rightarrow \Homo^n_{dR}(V^{an};\nabla_x)$ is defined by taking harmonic representatives.


By the uniqueness of harmonic representatives, $\Homo^n_{(2)}(V^{an},\mathcal{L}^\vee)$ is naturally endowed with a Hodge structure $\{\Homo^{p,q}\}_{p+q=n}$ of weight $n$. Let us take $[\omega],[\eta]\in \Homo^{n,0}\subset\Homo^n_{(2)}(V^{an},\mathcal{L}^\vee)$, that is, $[\omega]$ and $[\eta]$ are represented by $L^2$-harmonic $(n,0)$-forms. We are going to compute $\langle [\omega],[\bar{\eta}]\rangle$. Recall that any section of $\Db_{X^{an}}^{{\rm mod}D,(p,0)}(H)$ (resp. $\Db_{X^{an}}^{{\rm mod}D,(0,p)}(H)$) is harmonic if and only if it is a holomorphic (resp. anti-holomorphic) section. Therefore, for any algebraic $n$-differential forms $\omega,\eta\in\Omega^n(V^{alg})$, $\omega$ and $\bar{\eta}$ are both harmonic forms. This implies the equality $\iota[\bar{\eta}]=[\bar{\eta}]$. Setting ${\rm reg}:=can^{-1}$ and ${\rm reg}_{(2)}:=can_{(2)}^{-1}$, we have ${\rm reg}([\bar{\eta}])={\rm reg}_{(2)}([\bar{\eta}])$. We set $[\bar{\xi}]:={\rm reg}([\bar{\eta}])$. Since the diagram
\begin{equation}
\xymatrix{
\Homo^n_{(2)}(V^{an},\mathcal{L}^\vee)\otimes \Homo^n_c(V^{an},\mathcal{L}) \ar[d]_{{\rm id}\otimes |\Phi|^{-2}} \ar[r]^-{\langle\bullet,\bullet\rangle_{ch}}&\C \\
\Homo^n_{(2)}(V^{an},\mathcal{L}^\vee)\times \Homo^n_{(2)}(V^{an},\overline{\mathcal{L}^\vee}) \ar[ur]_-{\langle\bullet,\bullet\rangle} & 
}
\end{equation}
is commutative, we obtain identities $\langle\omega,|\Phi|^2\bar{\xi}\rangle_{ch}=\langle\omega,\bar{\xi}\rangle=\langle\omega,\bar{\eta}\rangle$.

On the other hand, we have an identification of local systems $\mathcal{L}^\vee\tilde{\rightarrow}\overline{\mathcal{L}}$ given by the correspondence $a\Phi^{-1}\mapsto a\overline{\Phi}$. With the aid of this identification, the homology intersection form is defined as a bilinear pairing
\begin{equation}
\langle\bullet,\bullet\rangle_{h}:\Homo_n\left( V^{an};\mathcal{L}\right)\times \Homo_n( V^{an};\overline{\mathcal{L}})\rightarrow\mathbb{C}.
\end{equation}

Combining the discussion above with the formula (\ref{eqn:expansion}), we obtain a
\begin{theorem}\label{thm:AbsoluteIntegral}
Let $\{[\Gamma_i]\}_{i=1}^r$ be a basis of $\Homo_n\left( V^{an};\mathcal{L}\right)$. One has a formula
\begin{equation}\label{AbsoluteIntegral}
I(\alpha)=\sum_{i,j}\left( \int_{\Gamma_i}\Phi\omega\right)C^{ji}\overline{\left( \int_{\Gamma_j}\Phi\eta\right)},
\end{equation}
where $(C^{ij})_{i,j}$ is the inverse of the intersection matrix $(\langle\Gamma_i,\overline{\Gamma_j}\rangle_h)$.
\end{theorem}

\begin{rem}
So far, we assumed that the parameters $\alpha_i$ are real. However, if we take into account the identity
\begin{equation}
\overline{\left( \int_{\Gamma_j}\Phi\eta\right)}= \int_{\overline{\Gamma_j}}\prod_{l=1}^k\overline{f_l(x)}^{\alpha_l}\overline{\eta} \ \ (\alpha_l\in\R),
\end{equation}
the right-hand side of (\ref{AbsoluteIntegral}) is clearly a meromorphic function in $\alpha_i\in\C$.
\end{rem}

\begin{example}
The simplest example of Theorem \ref{thm:AbsoluteIntegral} is when $A$ is given by $A=\begin{pmatrix}1&1\\ 0&1\end{pmatrix}$. The corresponding integral $\int_{\C}|t|^{2(\alpha-1)}|1-t|^{2(\beta-1)}dt\wedge d\bar{t}$ is well-known. We set $\Phi=t^\alpha (1-t)^\beta$ and $\omega=\eta=\frac{dt}{t(1-t)}$. Let $P\in\Homo_1(\C\setminus\{ 0,1\};\C\Phi)$ be the regularization (\cite[\S3.2]{AK}) of the interval $(0,1)$. If we set $e(\alpha)=e^{2\pi\ii\alpha}$, we obtain $\langle P,P^\vee\rangle_h=\frac{1-e(\alpha+\beta)}{(1-e(\alpha))(1-e(\beta))}$. Therefore, we have 
\begin{equation}
\int_{\C}|\Phi|^2\omega\wedge \bar{\omega}=\frac{(1-e(\alpha))(1-e(\beta))}{1-e(\alpha+\beta)}B(\alpha,\beta)^2.
\end{equation}
We write $t=\tau_1+\ii\tau_2$. Since $dt\wedge d\bar{t}=-2\ii d\tau_1\wedge d\tau_2$, we obtain
\begin{equation}\label{eqn:DoubleCopy}
\int_{\C}|t|^{2(\alpha-1)}|1-t|^{2(\beta-1)}d\tau_1\wedge d\tau_2=\frac{\sin\pi\alpha\sin\pi\beta}{\sin\pi(\alpha+\beta)}B(\alpha,\beta)^2.
\end{equation}
The formula (\ref{eqn:DoubleCopy}) was also discussed in \cite[(3.64)]{Dotsenko} and \cite[Corollary 1]{MimachiComplex}.
\end{example}

\begin{example}
We consider the complex Selberg integral discussed in \cite{AomotoComplexSelberg}. We set $n=N-2$ and consider an integral
\begin{equation}
I=\left( \frac{\ii}{2}\right)^{n}\int_{\C^{N-2}}\prod_{j=3}^N|z_j|^{2\alpha}|z_j-1|^{2\beta}\prod_{3\leq i<j\leq N}|z_i-z_j|^{2\gamma}dz_3\wedge\cdots\wedge dz_N\wedge d\bar{z}_3\wedge\cdots\wedge d\bar{z}_N
\end{equation}
We set $\Phi=\prod_{j=3}^Nz_j^{\alpha}(1-z_j)^{\beta}\prod_{3\leq i<j\leq N}(z_i-z_j)^{\gamma}$ and define an affine variety $V$ by $V:=\{z\in\C^n\mid \prod_{j=3}^Nz_j(z_j-1)\prod_{3\leq i<j\leq N}(z_i-z_j)\neq 0\}$. The symmetric group $\mathfrak{S}_n$ acts on $V$ by the permutation of the coordinates. Since $\s\nabla=\nabla\s$, $\mathfrak{S}_n$ also acts on the de Rham cohomology group $\Homo_{dR}^n(V,\nabla)$. We have
\begin{equation}
\dim_{\C}\Homo_{dR}^n(V,\nabla)=n!,\ \ \dim_{\C}\Homo_{dR}^n(V,\nabla)^{\mathfrak{S}_n}=1.
\end{equation}
The generator of the $\mathfrak{S}_n$-invariant part $\Homo_{dR}^n(V,\nabla)^{\mathfrak{S}_n}$ is given by the class $[\frac{dz_3\wedge\cdots\wedge dz_N}{\prod_{3\leq i<j\leq N}(z_i-z_j)}]$. Since $\mathfrak{S}_n$ defines a properly discontinuous action on $V$, the quotient morphism $\pi:V\rightarrow \mathfrak{S}_n\backslash V$ is a covering map and we obtain a canonical isomorphism
\begin{equation}\label{eqn:QuotIsom}
\Homo_{dR}^n(\mathfrak{S}_n\backslash V;\nabla)\simeq\Homo_{dR}^n(V;\nabla)^{\mathfrak{S}_n}
\end{equation}
induced by the pull-back $\pi^*$. Moreover, the basis of the $\mathfrak{S}_n$-invariant part of the twisted homology group $\Homo_n(V^{an};\mathcal{L})^{\mathfrak{S}_n}\simeq \Homo_n^{lf}(V^{an};\mathcal{L})^{\mathfrak{S}_n}$ is given by (the regularization of) a chamber $\Delta=\{ z\in\R^n\cap V^{an}\mid 0<z_i<1\ (i=3,\dots,N)\}$. Since $\Phi$ is $\mathfrak{S}_n$-invariant (up to a constant), the local system $\mathcal{L}$ induces a local system on $\mathfrak{S}_n\backslash V^{an}$ which is also denoted by $\mathcal{L}$ by abuse of notation. The dual of the isomorphism (\ref{eqn:QuotIsom}) is given by
\begin{equation}\label{eqn:QuotIsom2}
\Homo_n(V^{an};\mathcal{L})\simeq\Homo_n(\mathfrak{S}_n\backslash V^{an};\mathcal{L})^{\mathfrak{S}_n}
\end{equation}
induced by the push-forward $\pi_*$. Therefore, Theorem \ref{thm:AbsoluteIntegral} applied to the de Rham cohomology group $\Homo_{dR}^n(\mathfrak{S}_n\backslash V;\nabla)\simeq\Homo_{dR}^n(V;\nabla)^{\mathfrak{S}_n}$ and combined with the result of \cite[Theorem 1]{MY} gives a formula
\begin{equation}
I=\frac{\prod_{j=1}^n\sin\pi(\alpha+\frac{(j-1)}{2}\gamma)\sin\pi(\beta+\frac{(j-1)}{2}\gamma)\sin\pi(\alpha+\frac{j\gamma}{2})}{n!\prod_{j=1}^n\sin\pi(\alpha+\beta+\frac{(n+j-2)}{2}\gamma)\sin\pi(\frac{\gamma}{2})}S(\alpha,\beta,\gamma)^2.
\end{equation}
Here, $S(\alpha,\beta,\gamma)$ is the ordinary Selberg integral 
\begin{align}
S(\alpha,\beta,\gamma)&=\int_{[0,1]^n}\prod_{j=3}^Nx_j^{\alpha}(x_j-1)^{\beta}\prod_{3\leq i<j\leq N}(x_i-x_j)^{\gamma}dx_3\wedge\cdots\wedge dx_N\\
&=\frac{\prod_{j=1}^n\Gamma(\alpha+1+\frac{(j-1)}{2}\gamma)\Gamma(\beta+1+\frac{(j-1)}{2}\gamma)\Gamma(\frac{j\gamma}{2}+1)}{\prod_{j=1}^n\Gamma(\alpha++\beta+2+\frac{(n+j-2)}{2}\gamma)\Gamma(\frac{\gamma}{2}+1)}.
\end{align}
This result is in concordance with the main result of \cite{AomotoComplexSelberg}.

\end{example}

\subsection{GKZ case}\label{sec:62}
We use the same notation as \S\ref{sec:GKZ}. Let us fix a unimodular regular triangulation $T$. For any $(n+k)$-simplex $\s$, we set
\begin{equation}\label{seriespsi}
\psi_{\s}(z;\delta)=
\sum_{{\bf m}\in\Z_{\geq 0}^{\bs}}\frac{z_\sigma^{-A_\sigma^{-1}A_{\bar{\sigma}}{\bf m}}z_{\bar{\sigma}}^{\bf m}}{\Gamma({\bf 1}_\sigma-A_\sigma^{-1}(\delta+A_{\bar{\sigma}}{\bf m})){\bf m!}}.
\end{equation}
In view of the formula (\ref{eqn:SigmaIntersection}), we obtain a

\begin{theorem}\label{thm:AbsoluteCIN}
For a unimodular regular triangulation $T$\footnote{This assumption can be removed by modifying the right-hand side of the formula (\ref{eqn:AbsoluteCIN})}, we have an identity
\begin{align}
&\left(\frac{\ii}{2}\right)^n\int_{\C^n}\prod_{l=1}^k|h_l(x;z)|^{-2\gamma_l}\prod_{i=1}^n|x_i|^{2c_i}\frac{dx}{x}\wedge\frac{d\bar{x}}{\bar{x}}\nonumber\\
=&\Gamma({\bf 1}-\gamma)^2\sin\pi\gamma\sum_{\s\in T}\frac{\pi^{2n}}{\sin\pi A_\s^{-1}\delta}
|z_\s|^{-2A_\s^{-1}\delta}
\psi_{\s}\left( z;\delta\right)
\psi_{\s}\left( \bar{z};\delta\right).\label{eqn:AbsoluteCIN}
\end{align}
The right-hand side of (\ref{eqn:AbsoluteCIN}) is convergent for any $z$ such that $(-\log|z_1|,\dots,-\log|z_N|)$ is in a sufficiently far translation of the cone $C_T$ inside itself.
\end{theorem}

\noindent
Note that $\psi_{\s,{\bf k}}\left( z;\delta\right)
\psi_{\s,{\bf k}}\left( \bar{z};\delta\right)=|\psi_{\s,{\bf k}}\left( z;\delta\right)|^2$ when $\delta$ is real.

\begin{example}
The second simplest example is Gau\ss' ${}_2F_1$ case. We set 
\begin{equation}
A=
\begin{pmatrix}
1&0&1&0\\
0&1&0&1\\
0&0&1&1
\end{pmatrix}.
\end{equation}
Taking the regular triangulation $T=\{ 124,134\}$ and substituting $z_1=1,z_2=z,z_3=z_4=-1$, we obtain
\begin{align}
&\frac{\ii}{2}\int_\C |1-x|^{-2\gamma_1}|z-x|^{-2\gamma_2}|x|^{2(c-1)}dx\wedge d\bar{x}\nonumber\\
=&\frac{\Gamma(1-\gamma_1)^2\Gamma(1-\gamma_2)^2}{\pi^2}\left\{ \frac{\sin^2\pi\gamma_1\sin\pi\gamma_2\sin\pi c}{\sin\pi(\gamma_2-c)}|z|^{2(c-\gamma_2)}{}_2f_1\left(\substack{\gamma_1,c\\ 1-\gamma_2+c};z\right){}_2f_1\left(\substack{\gamma_1,c\\ 1-\gamma_2+c};\bar{z}\right)\right.\nonumber\\
&\left.+\frac{\sin\pi\gamma_1\sin^2\pi\gamma_2\sin\pi (\gamma_1+\gamma_2-c)}{\sin\pi(c-\gamma_2)}
{}_2f_1\left(\substack{\gamma_1+\gamma_2-c,\gamma_2\\ 1-c+\gamma_2};z\right)
{}_2f_1\left(\substack{\gamma_1+\gamma_2-c,\gamma_2\\ 1-c+\gamma_2};\bar{z}\right)\right\},
\end{align}
where ${}_2f_1\left(\substack{\alpha,\beta\\ \gamma};z\right)$ is given by the formula
\begin{equation}
{}_2f_1\left(\substack{\alpha,\beta\\ \gamma};z\right)
=
\sum_{m=0}^\infty\frac{\Gamma(\alpha+m)\Gamma(\beta+m)}{\Gamma(\gamma+m)m!}z^m. 
\end{equation}
This formula is equivalent to \cite[Corollary 2]{MimachiComplex}. Indeed, the relation between our parameter $\delta=\begin{pmatrix}\gamma_1\\ \gamma_2\\ c\end{pmatrix}$ and the parameters $a,b,c$ of \cite[Corollary 2]{MimachiComplex} is given by
\begin{equation}
\delta=
\begin{pmatrix}
-b\\
-c\\
a+1
\end{pmatrix}.
\end{equation}
The same formula with different notation is also obtained in \cite[(3.64)]{Dotsenko}.
\end{example}

\section*{Acknowledgement}

This article is a detailed version of a talk presented at the workshop ``MathemAmplitudes 2019: Intersection Theory and Feynman Integrals" held in Padova, Italy on 18-20 December 2019. The author thanks the participants and organizers of this meeting. The results of \S6 is an outcome of a discussion with Sebastian Mizera. The author would like to thank him. The author thanks his collaborators Yoshiaki Goto and Nobuki Takayama for many discussions and valuable comments.



\end{document}